\documentclass[review]{siamonline220329}
\usepackage {epsfig,amsmath,euscript}
\usepackage[latin1]{inputenc}
\usepackage[english]{babel}
\usepackage{color}
\usepackage{amsmath}
\usepackage{amssymb}
\usepackage{ae}
\usepackage{subfigure}
\usepackage[T1]{fontenc}
\usepackage{graphicx}
\usepackage{geometry}
\usepackage{epstopdf}
\usepackage{multicol}
\usepackage{hyperref}
\usepackage{color}
\usepackage{makecell}
\usepackage{tikz}
\usepackage{circuitikz}
\usepackage{overpic}

\definecolor{rred}{rgb}{0.7,0.0,0.2}
\definecolor{bblue}{rgb}{0.2,0.0,0.7}
\usepackage{tikz}
\usepackage{rotating}
\usepackage{todonotes}
\usepackage{bm}
\usepackage{colortbl}
\usepackage{tabularx}
\usepackage{xcolor}
\newcommand*\circled[1]{\tikz[baseline=(char.base)]{
            \node[shape=circle,draw,inner sep=2pt] (char) {#1};}}

\usepackage{comment} 

\title{The effects of delay on the HKB model of human motor coordination\thanks{Submitted to the editors 28 October 2022.
\funding{SJH would like to thank the Hungarian Academy of Sciences for support through its Distinguished Guest Scientist Programme.}}}

\nolinenumbers
\author{L. I. Allen\thanks{Population Health Sciences, University of Bristol, United Kingdom (\email{hk19479@bristol.ac.uk}).} 
\and T. G. Moln\'ar\thanks{Department of Mechanical and Civil Engineering, California Institute of Technology, Pasadena, CA 91125, USA (\email{tmolnar@caltech.edu}).}
\and Z. Domb\'ov\'ari\thanks{MTA-BME Lend{\"{u}}let Machine Tool Vibration Research Group, Department of Applied Mechanics, Faculty of Mechanical Engineering, Budapest University of Technology and Economics, Budapest 1111, Hungary (\email{dombovari@mm.bme.hu}).}
\and S. J. Hogan\thanks{Department of Engineering Mathematics, University of Bristol, Bristol BS8 1UB, United Kingdom (\email{s.j.hogan@bristol.ac.uk}). Corresponding Author: ORCiD:  0000-0001-6012-6527}}

\headers{The effects of delay on the HKB model of human motor coordination}{L. I. Allen, T. G. Molnar, Z. Dombovari, and S. J. Hogan}

\date{Received: date / Accepted: date}

\begin{document}
\maketitle
\nolinenumbers
\begin{abstract}
In this paper, we analyse the celebrated Haken-Kelso-Bunz (HKB) model, describing the dynamics of bimanual coordination, in the presence of delay.
We study the linear dynamics, stability, nonlinear behaviour and bifurcations of this model by both theoretical and numerical analysis.
We calculate in-phase and anti-phase limit cycles as well as quasi-periodic solutions via double Hopf bifurcation analysis and centre manifold reduction.
Moreover, we uncover further details on the global dynamic behaviour by numerical continuation, including the occurrence of limit cycles in phase quadrature and 1-1 locking of quasi-periodic solutions.

\end{abstract}



\section{Introduction}
Human motor coordination is the result of complex interactions, at many different time and length scales. One established way to model bimanual coordination is to assume that the fingers or limbs of experimental subjects are oscillators, capable of generating self-sustained periodic motion.
Within this approach, the focus is on understanding the observed {\it relative phase} $\phi$ of the two oscillators \cite{Kelso1995}. 
Stable \textit{in-phase} ($\phi=0$) synchronisation is usually found to be the simplest to maintain \cite{Buchanan2006}, \cite{Kelso1981}, while stable \textit{anti-phase} ($\phi=\pi$) motion  \cite{Bourbousson2010} and stable \textit{phase-lagged} states can also occur \cite{Collins1993}, \cite{Duarte2012}.

The Haken-Kelso-Bunz (HKB) coupled oscillator model \cite{Haken1985, Kelso2021} was originally developed to explain these different types of bimanual synchronisation, and the ways in which they can occur. Since then it has become the bedrock of all subsequent research in this area. 

In more recent times, the HKB model has found application elsewhere. Schizophrenia is a mental illness with high prevalence and a scarcity of satisfactory treatments. This has inspired a  drive within the research community \cite{Slowinski2017} to develop methods for early diagnosis and preventative intervention.

The mirror game, where two individuals mirror each other's movements, is considered to be a powerful tool for studying coordination dynamics \cite{Noy2011}. In experiments based on this principle, Varlet et al.~\cite{Varlet2012} identified peculiar characteristics of the motion of schizophrenic patients performing simple synchronisation tasks with healthy individuals. 
In one version of the mirror game, an HKB-driven virtual player \cite{Zhai2018} participates with a human partner. The analysis of the resulting motion could form the basis of a diagnostic tool for schizophrenia \cite{Lombardi2018}, where motor abnormalities are one of the first indicators of the illness. 

In their work, Varlet et al.~\cite{Varlet2012} pointed out that the standard HKB coupled oscillator model \cite{Haken1985} does not accurately describe the observed dynamics of schizophrenia patients unless a delay is included in the coupling term. The resulting {\it delayed HKB model} \cite{Banerjee2006, Slowinski2016, Washburn2015, Zhai2021} is the subject of this paper. Our goal is to inform further research into human coordination dynamics where delay is especially prevalent.
In particular, we seek to aid the understanding of the delayed HKB equation by providing insight into the underlying dynamic behaviour of the system. Varlet et al.~\cite{Varlet2012} suggest that the delayed HKB equation may be suitable to capture the mirror game whereas the numerical approaches used by S\l{}owi\'nski et al.~\cite{Slowinski2020} indicate that the delayed HKB equation may have limited relevance to experimentally observed behaviour. We aim to help clarify the extent of the utility of the delayed HKB model\footnote{We note that in the analysis by S\l{}owi\'nski et al., various parameters are fixed at values found by Kay et al.~\cite{Kay1987} using experiments studying the hand motion of four participants. However, Peper et al.~\cite{Peper2004} suggest that not all limbs can be modelled in the same way which motivates a broader investigation of the parameter space, as presented here.}. 

Our paper is organised as follows.
Section~\ref{chap:HKBeqns} describes the delayed HKB model of human motor coordination.
Section~\ref{chap:linear_stability_charts} presents linear stability analysis that is verified by numerical results in Section~\ref{sec:stab_chart_comp}.
The nonlinear dynamics of the delayed HKB model are analysed, by centre manifold reduction, in Section~\ref{chap:bifurcation_analysis}.
We compare these theoretical results with numerical continuation
in Section~\ref{chap:numerics}.

\section{The delayed HKB model}
\label{chap:HKBeqns}

The delayed HKB model \cite{Banerjee2006, Slowinski2016} is given by

\begin{equation}
\label{HKBeqns:full_delayed_HKB}
\begin{aligned}
\ddot{x}_1(t)+\omega^2x_1(t)&=\left(\gamma-\alpha x_1^2(t)-\beta \dot{x}_1^2(t)\right)\dot{x}_1(t)\\
& \quad +\left(a+b\left(x_1(t)-x_2\left(t-\tau_1\right)\right)^2\right)
\left(\dot{x}_1(t)-\dot{x}_2\left(t-\tau_1\right)\right), \\
\ddot{x}_2(t)+\omega^2x_2(t)&=\left(\gamma-\alpha x_2^2(t)-\beta \dot{x}_2^2(t)\right)\dot{x}_2(t)\\
& \quad +\left(a+b\left(x_2(t)-x_1\left(t-\tau_2\right)\right)^2\right)
\left(\dot{x}_2(t)-\dot{x}_1\left(t-\tau_2\right)\right).
\end{aligned}
\end{equation}
This model is a pair of coupled second-order delay differential equations (DDEs). The variables $x_1(t)$ and $x_2(t)$ represent the amplitudes\footnote{In the original HKB paper \cite{Kelso1995}, $x_1(t)$ and $x_2(t)$ are the angular displacements of each finger, with direction defined symmetrically so that $x_1(t)=x_2(t)$ corresponds to in-phase motion.} of the individual oscillators at time $t$. The parameter $\gamma$ is the linear damping coefficient and $\alpha$, $\beta$ are nonlinear damping coefficients, also known as the Van der Pol and Rayleigh coefficients, respectively. Parameter $a$ is the linear coupling coefficient and $b$ is the nonlinear coupling coefficient. The pacing frequency $\omega$ is physically positive\footnote{In \cite{Kelso1995}, pacing was provided by a metronome.}.
The time delays $\tau_1$, $\tau_2$ arise from cognitive and physiological processes, typically caused by detection and actuation, which can be different for each oscillator. In this paper, we take these two time delays to be equal: $\tau_1=\tau_2=:\tau$. 
We remark that the number of parameters in the delayed HKB model \eqref{HKBeqns:full_delayed_HKB} could be reduced by introducing the scaled time ${\tilde{t} = \omega t}$, however, hereinafter we rather use functions of $t$ for easier physical interpretation.

The delayed HKB model \eqref{HKBeqns:full_delayed_HKB} has discrete symmetries in its structure: $x_1$ and $x_2$, as well as $x_1$ and $-x_2$ are interchangeable.
These symmetries will ultimately result in the existence of in-phase and anti-phase periodic solutions.

We begin our analysis of the delayed HKB model \eqref{HKBeqns:full_delayed_HKB} by first considering the linearized version of the system, given\footnote{To simplify notation, we do not denote the dependence of $x_1$ and $x_2$ on time, unless the delay is involved.} by
\begin{equation}
\label{HKBeqns:linear_delay_HKB}
\begin{aligned}
\ddot{x}_1+\omega^2x_1&=\gamma\dot{x}_1+a\left(\dot{x}_1-\dot{x}_2(t-\tau)\right), \\
\ddot{x}_2+\omega^2x_2&=\gamma\dot{x}_2+a\left(\dot{x}_2-\dot{x}_1(t-\tau)\right).
\end{aligned}
\end{equation}
We shall show that these equations provide us with an explanation of the fundamental structures observed in numerical computations \cite{Slowinski2016}.

In \cite[eq. (5)]{Cass2021}, it was shown that the linear HKB equations in the \textit{absence} of delay could be simplified when written in terms of normal modes. We adopt the same approach here, by setting $\eta^{(i)}=x_1+x_2$ and $\eta^{(a)}=x_1-x_2$, corresponding to in-phase motion and anti-phase motion respectively, so that \eqref{HKBeqns:linear_delay_HKB} becomes

\begin{equation}
\label{HKBeqns:normal_forms}
\begin{aligned}
\ddot{\eta}^{(i)}+\omega^2\eta^{(i)}&=\gamma\dot{\eta}^{(i)}+a\left(\dot{\eta}^{(i)}-\dot{\eta}^{(i)}(t-\tau)\right), \\
\ddot{\eta}^{(a)}+\omega^2\eta^{(a)}&=\gamma\dot{\eta}^{(a)}+a\left(\dot{\eta}^{(a)}+\dot{\eta}^{(a)}(t-\tau)\right).
\end{aligned}
\end{equation}

We consider the stability of the trivial solutions (equilibria) $\eta^{(i,a)}=0$ of \eqref{HKBeqns:normal_forms}. 
\begin{itemize}
    \item If both $\eta^{(i,a)}=0$ are stable, then the equilibrium $x_1=x_2=0$ of \eqref{HKBeqns:linear_delay_HKB} must be stable and we will see no oscillations. 
    \item When $\eta^{(i)}=0$ is unstable and $\eta^{(a)}=0$ is stable, we expect to find stable in-phase limit cycles in the full system \eqref{HKBeqns:full_delayed_HKB}. 
    \item When $\eta^{(i)}=0$ is stable and $\eta^{(a)}=0$ is unstable, we expect to find stable anti-phase limit cycles in the full system \eqref{HKBeqns:full_delayed_HKB}. 
\end{itemize}

Such limit cycles arise because of a Hopf bifurcation, which occurs when changing parameters of the system.
If the leading eigenvalues $\lambda$ cross the imaginary axis in pairs, $\lambda=\pm {\rm i}\nu$ where $\nu\neq0$, then, according to an infinite dimensional version of the Hopf Bifurcation Theorem \cite{Hassard1981}, there will be a Hopf bifurcation when the derivative of the eigenvalue $\lambda$ with respect to a parameter evaluated at $\lambda={\rm i}\nu$ is non-zero.

Equations \eqref{HKBeqns:normal_forms} are uncoupled. Since they only differ by one sign, we can carry out an analysis on the \textit{linear delayed HKB model}, given by
\begin{equation}
\label{HKB:normal_forms_together}
\ddot{\eta}+\omega^2\eta=\gamma\dot{\eta}+a\left(\dot{\eta}\mp\dot{\eta}(t-\tau)\right),
\end{equation}
where $\eta = \eta^{(i)}$ corresponds to the minus sign in the last term and $\eta = \eta^{(a)}$ corresponds to the plus sign. 

\section{Stability charts for the linear delayed HKB model}
\label{chap:linear_stability_charts}

In this section we produce stability charts in the $(a,\tau)$, $(\gamma,\tau)$ and $(\gamma,a)$ planes, for arbitrary values of $\omega$, of the trivial solutions of \eqref{HKB:normal_forms_together}. These charts will also be used to explain structure of the bifurcation diagram for the nonlinear delayed HKB model \eqref{HKBeqns:full_delayed_HKB} in the $(a,\tau)$ plane produced by S\l{}owi\'{n}ski et al. \cite[Figure 2(a)]{Slowinski2016}.

\begin{figure}
	\centering
	\includegraphics[width=.9\columnwidth]{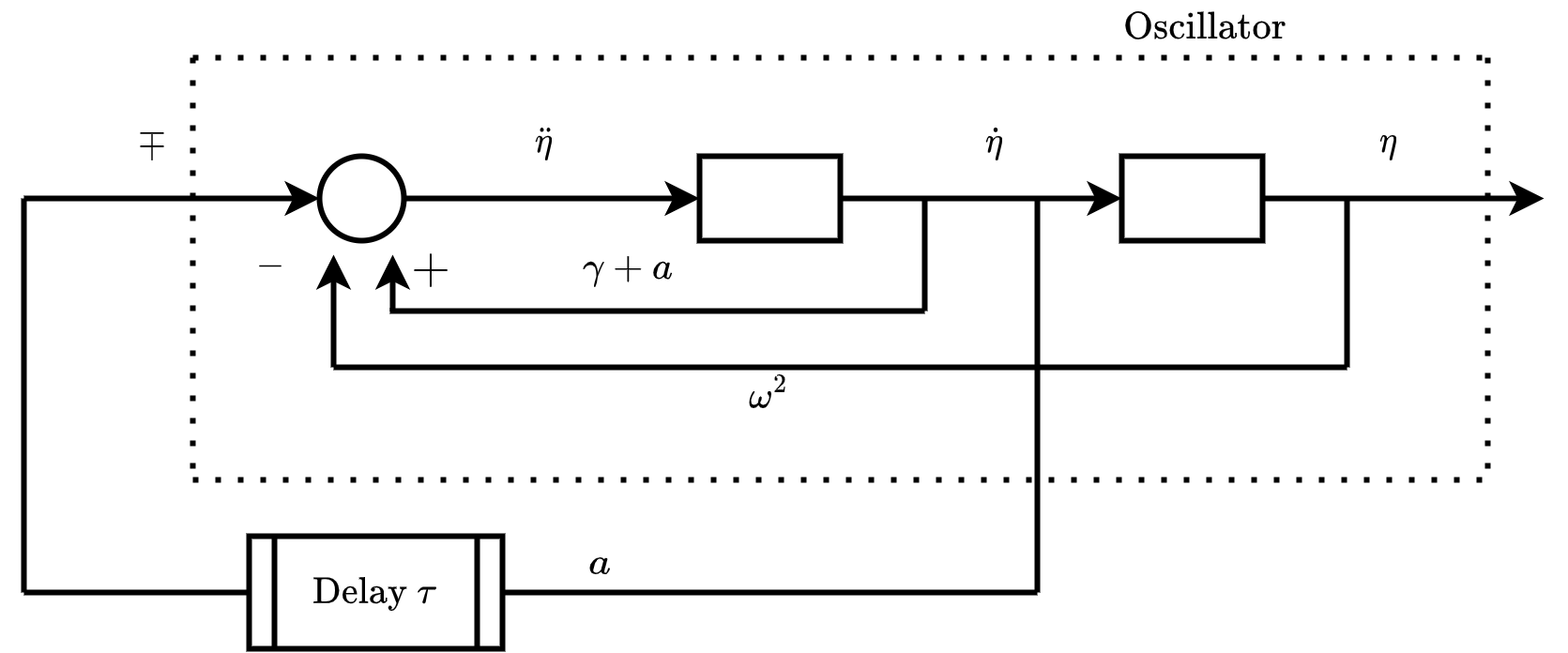}
	\caption[Mechanical analogue for the normal modes of the delayed HKB model]{A feedback control system block diagram for the linear delayed HKB model \eqref{HKB:normal_forms_together_control_form}. 
	The damped harmonic oscillator on the left-hand side of \eqref{HKB:normal_forms_together_control_form} is within the dotted section.} 
	\label{fig:intro:mechanicalanalogue}
\end{figure}

Equations \eqref{HKB:normal_forms_together} will be very familiar to control engineers when written in the form:
\begin{equation}
\label{HKB:normal_forms_together_control_form}
\ddot{\eta}(t)-(\gamma+a)\dot{\eta}(t)+\omega^2\eta(t)=\mp a \dot{\eta}(t-\tau),
\end{equation}
The left-hand side of \eqref{HKB:normal_forms_together_control_form} represents damped simple harmonic motion, with damping coefficient\footnote{We shall consider cases when $-(\gamma+a) \gtrless 0$. } $-(\gamma+a)$. The right-hand side is a delayed harmonic oscillator with feedback gain $\mp a$, which changes the damping coefficient $-(\gamma+a)$. The feedback control system block diagram  is shown in Figure~\ref{fig:intro:mechanicalanalogue}.

In the absence of delay, the in-phase and anti-phase trivial solutions $\eta^{(i,a)}=0$ have two lines of Hopf-bifurcations in $(a,\gamma)$ parameter space; $HB_I: \gamma=0$ and $HB_A:2a+\gamma=0$, see \cite[Figure 1]{Cass2021}.  Bistable regions of these normal modes are also seen.

In the presence of delay, a stability analysis of the equilibrium solution $\eta=0$ to \eqref{HKB:normal_forms_together_control_form} is more complicated, but well established \cite{Bhatt1966, Hsu1966, Kuang1993}. Details are given in Appendix~\ref{appA}.

\begin{figure}
	\centering
	\includegraphics[scale=1]{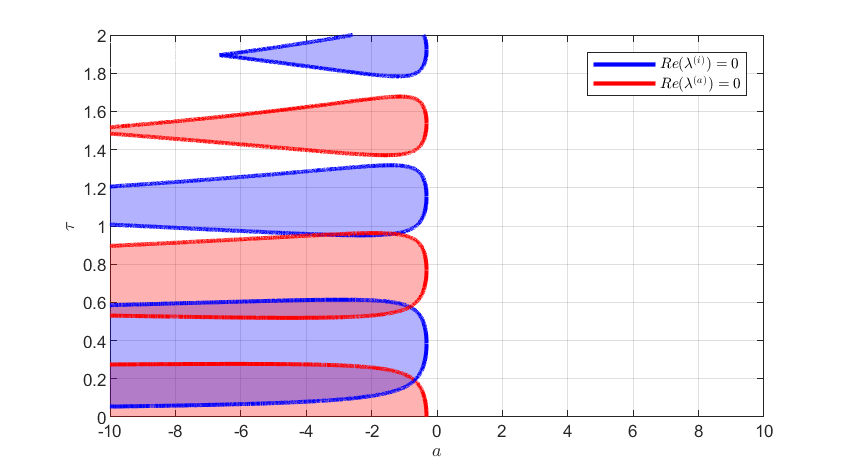}
	\caption[$(a,\tau)$ stability chart for the normal modes of the delayed HKB system]{Stability chart in $(a,\tau)$ parameter space for the normal modes of the linear delayed HKB equation \eqref{HKB:normal_forms_together_control_form} with $\gamma=0.641,\ \omega=2.6\pi$. Regions where $\eta^{(i)}$ is stable are shaded blue, regions where $\eta^{(a)}$ is stable are shaded red. The normal modes $\eta^{(i,a)}$ are both stable in the purple regions and both unstable in the white regions.}
	\label{fig:linear_stability_charts:Kuang_a_tau}
\end{figure}

For the stability chart in $(a, \tau)$ parameter space, shown in Figure~ \ref{fig:linear_stability_charts:Kuang_a_tau}, we follow \cite{Slowinski2016} and set\footnote{We reproduced their results by setting $\omega=2\pi(1.3)=2.6\pi$, suggesting that they took $\omega=1.3[Hz]$.} $\gamma=0.641$, $\omega=2.6\pi$. Regions where $\eta^{(i)}=0$ is stable are shaded blue, regions where $\eta^{(a)}=0$ is stable are shaded red. Both $\eta^{(i,a)}=0$ are stable in the purple regions. The eigenvalues of the in-phase and anti-phase normal modes are denoted by $\lambda^{(i,a)}$ respectively. Stability boundaries occur when the rightmost eigenvalues are pure imaginary: $\lambda^{(i,a)} = {\rm i}\nu^{(i,a)}$. 
From \eqref{rho_definition}, the sign of the quantity ${\rho_{\tau}=\text{Re}\left(\frac{\partial\lambda}{\partial\tau}\bigr\rvert_{\lambda ={\rm i}\nu}\right)}$ indicates how the stability of $\eta=0$ changes as $\tau$ increases. 

In the absence of delay, along the line $\tau=0$ in Figure~ \ref{fig:linear_stability_charts:Kuang_a_tau}, $\eta^{(a)}=0$ is stable (and $\eta^{(i)}=0$ is unstable) for $a <-\gamma/2= -0.3205$, in agreement with the stability boundary $HB_A:2a+\gamma=0$ in \cite[Figure 1]{Cass2021}. We observe finite amplitude in-phase limit cycles in the full equations \eqref{HKBeqns:full_delayed_HKB} with $\tau=0$ \cite{Avitabile2016, Cass2021}. 

As $\tau$ increases for fixed $a <-\gamma/2= -0.3205$, initially $\eta^{(a)}=0$ remains stable. Then we observe a number of stability switches. $\eta^{(i)}=0$ becomes stable on crossing the lowest blue line in Figure~ \ref{fig:linear_stability_charts:Kuang_a_tau}. When both $\eta^{(i,a)}=0$ are stable (in the purple regions), we have the case when $x_1=x_2=0$ is stable. So, small but finite values of the delay $\tau$ eliminate the finite amplitude in-phase limit cycles in \eqref{HKBeqns:full_delayed_HKB}. On crossing the lowest red line in Figure~ \ref{fig:linear_stability_charts:Kuang_a_tau}, $\eta^{(a)}=0$ loses stability. We expect to see finite amplitude anti-phase limit cycles in the full delayed HKB model \eqref{HKBeqns:full_delayed_HKB} in the blue region around $\tau=0.4$ for $a <-\gamma/2= -0.3205$. Similar observations can be made as we increase $\tau$ further in Figure~ \ref{fig:linear_stability_charts:Kuang_a_tau}.

\begin{figure}[ht!]
	\centering
	\includegraphics[width=\columnwidth]{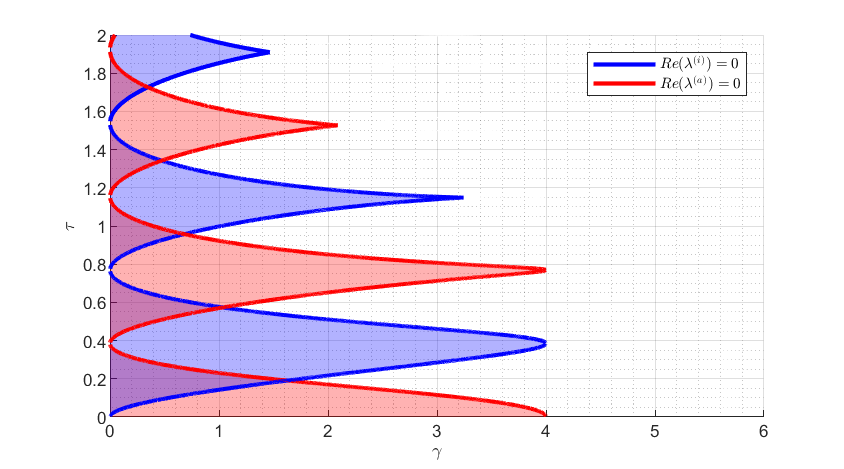}
	\caption[$(\gamma,\tau)$ stability chart for the normal modes of the delayed HKB system]{Stability chart in $(\gamma,\tau)$ parameter space for the normal modes of the linear delayed HKB equation \eqref{HKB:normal_forms_together_control_form} with $a=-2,\ \omega=2.6\pi$. 
    The colour scheme is the same as that in Figure~\ref{fig:linear_stability_charts:Kuang_a_tau}.
	}
	\label{fig:linear_stability_charts:Kuang_gamma_tau}
	\vspace{-5mm}
\end{figure}

For the stability chart in $(\gamma,\tau)$ parameter space, shown in Figure~
\ref{fig:linear_stability_charts:Kuang_gamma_tau}, we set $a=-2$ and $\omega=2.6\pi$ \cite{Slowinski2016}. The colour scheme and the definitions of $\lambda^{(i,a)}$ and $\rho^{(i,a)}_\tau$ are the same as those in Figure~\ref{fig:linear_stability_charts:Kuang_a_tau}. In Figure~\ref{fig:linear_stability_charts:Kuang_gamma_tau}, along the line $\tau=0$, we see that $\eta^{(a)}=0$ is stable (and $\eta^{(i)}=0$ is unstable) for $\gamma < -2a= 4$, in agreement with the stability boundary $HB_A:2a+\gamma=0$ in \cite[Figure 1]{Cass2021}. We see further stability switches as $\tau$ increases. 

Figures~\ref{fig:linear_stability_charts:Kuang_a_tau} and \ref{fig:linear_stability_charts:Kuang_gamma_tau} show a periodic nature in the stability curves as $\tau$ increases. 
We can explain this observation as follows. In \eqref{linear_stability_charts:HKBlin_Dcurve_special_case}, we show  that when $\tau = \frac{n\pi}{\omega}$, $n\in\mathbb{Z}$, stability boundaries are given by $\gamma=a(\pm(-1)^n-1)$ for $\omega\neq0$. When $n$ is even, that means the boundaries are $\gamma=0$ for $\eta^{(i)}$ and $\gamma+2a=0$ for $\eta^{(a)}$. When $n$ is odd, the boundaries are swapped and given by $\gamma+2a=0$ for $\eta^{(i)}$ and $\gamma=0$ for $\eta^{(a)}$. 

\begin{figure}[]
\centering
\begin{subfigure}
\centering
\includegraphics[width=.45\textwidth]{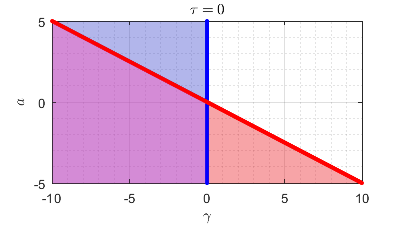}
\label{fig:gamma_a:fig_0}
\end{subfigure}
\begin{subfigure}
\centering
\includegraphics[width=.45\textwidth]{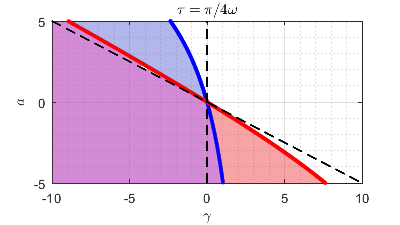}
\label{fig:gamma_a:fig_1_a}
\end{subfigure}%
\begin{subfigure}
\centering
\includegraphics[width=.45\textwidth]{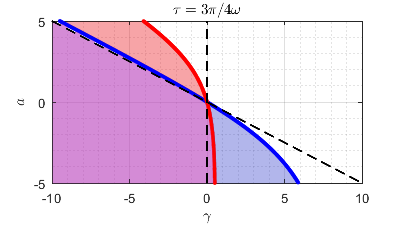}
\label{fig:gamma_a:fig_1_c}
\end{subfigure}
\begin{subfigure}
\centering
\includegraphics[width=.45\textwidth]{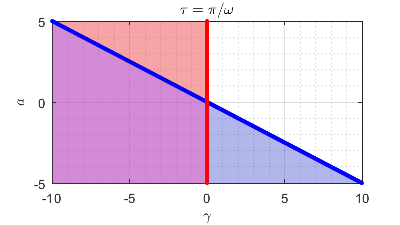}
\label{fig:gamma_a:fig_1_d}
\end{subfigure}
\begin{subfigure}
    \centering
    \includegraphics[scale = 1]{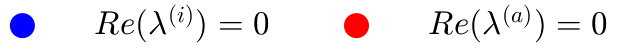}
\end{subfigure}%
\caption[$(\gamma, a)$ stability charts for the normal modes of the delayed HKB system, small $\tau$]{
	Stability chart in $(\gamma,a)$ parameter space for the normal modes of the linear delayed HKB equation \eqref{HKB:normal_forms_together_control_form}, for different values of $\tau$, with $\omega=2.6\pi$.
	The colour scheme is the same as that in Figure~\ref{fig:linear_stability_charts:Kuang_a_tau}.
	The stability chart with $\tau=0$ was given in \cite[Figure 1]{Cass2021}, where regions of {\it instability} were highlighted. Within the parameter ranges illustrated, the stability boundaries for $\tau=\frac{\pi}{\omega}$ are the same as those for $\tau=0$, but with the in-phase and anti-phase curves swapped.}
	\label{fig:gamma_a:fig1}
\end{figure}

Stability charts in the $(\gamma,a)$ plane are illustrated in Figure~\ref{fig:gamma_a:fig1}, 
showing how the stability boundaries evolve as $\tau$ increases from $0$ to $\pi/\omega$. At $\tau=0$ the in-phase and anti-phase stability curves are $\gamma=0$ and $2a+\gamma=0$ respectively. As $\tau$ increases, these boundaries deform and cross, eventually switching when $\tau=\pi/\omega$, as expected.

\section{Comparison with numerical stability charts}
\label{sec:stab_chart_comp}

Up to now, we have analysed the linear delayed HKB  \eqref{HKBeqns:linear_delay_HKB}, which is valid for small amplitudes. However, most experiments are performed in the nonlinear regime. The study by S\l{}owi\'{n}ski et al. \cite{Slowinski2016} is the only bifurcation analysis of the nonlinear delayed HKB model \eqref{HKBeqns:full_delayed_HKB}, but it is entirely numerical, and does not consider the linearised equations. These authors fixed the linear damping coefficient $\gamma=0.641$ and the nonlinear damping coefficients $\alpha=12.457$, $\beta=0.007905$ \cite{Kay1987}. Variation of the other parameters was considered within the range of experimentally observed values. We follow these authors with the same choice of parameters, with  $\omega=2\pi(1.3)=2.6\pi$.

Figure~\ref{fig:linear_stability_charts:Kuang_gamma_tau} shows the analytic stability curves in $(\gamma,\tau)$ parameter space. Our own numerical continuation of the Hopf bifurcations using \textsc{DDE-Biftool} \cite{DDEBifTool2002, DDEBifToolMan2017} shows exact agreement between the numerical and analytic curves (not shown). Since the analytic curves lie exactly on the numerical curves, the generic Hopf bifurcation curves can be deduced exactly from analysis of the linear system \eqref{HKBeqns:linear_delay_HKB}.

Figure~\ref{fig:linear_stability_charts:Kuang_gamma_tau} shows that the stability curves satisfy $0<\gamma<4$. This can be explained using the analysis in Appendix~\ref{appA}. Conditions \eqref{linear_stability_charts:conditionImRootmapped} for eigenvalues to touch or cross the imaginary axis can be combined to give $\gamma(\gamma+2a)<0$. Therefore, for $a=-2$ as in Figure~\ref{fig:linear_stability_charts:Kuang_gamma_tau}, we have $0<\gamma<4$.

S\l{}owi\'{n}ski et al. \cite[Figure 2(a)]{Slowinski2016} carried out a comprehensive numerical analysis of the full problem \eqref{HKBeqns:full_delayed_HKB}, using \textsc{DDE-Biftool}. They found many bifurcations in\footnote{Experiments tend to suggest that $a<0$.} $a>0$, and torus bifurcations in $a<0$.

\section{Applying centre manifold theory to the delayed HKB equation}
\label{chap:bifurcation_analysis}

A standout observation of Figures~\ref{fig:linear_stability_charts:Kuang_a_tau} and \ref{fig:linear_stability_charts:Kuang_gamma_tau} is the crossing of Hopf bifurcation curves. At such points we expect to find double Hopf (or Hopf-Hopf) bifurcations, as two pairs of eigenvalues cross the imaginary axis at the same time.

We now revert to the full problem \eqref{HKBeqns:full_delayed_HKB}, and investigate the double Hopf points in $(a,\tau)$ parameter space, unfolding the in-phase and anti-phase periodic orbits as parameter values are varied nearby. 
We use centre manifold theory, with symbolic computations in Maple\textsuperscript{TM} based on the tutorial codes given in~\cite{Campbell2009}; the details of this calculation are outlined in Appendix~\ref{appB}.


Using \eqref{linear_stability_charts:Kuang_tau} and \eqref{linear_stability_charts:Kuang_sincos},  we find four double Hopf bifurcations at $(a,\tau)=(a_c,\tau_c)$ for\footnote{This range was given in Figure~\ref{fig:linear_stability_charts:Kuang_a_tau}.} $a\in[-10,10]$, $\tau\in(0,2]$, labelled  \texttt{HH1} to \texttt{HH4} in Table \ref{tab:Hopf_points}, along with values for the critical eigenvalues $\lambda^{(i,a)}={\rm i}\nu^{(i,a)}$.  

\begin{table*}[ht]
	\centering
	\def\arraystretch{1.5}
	\setlength{\arrayrulewidth}{1.5pt}
	\small
\begin{tabularx}{\textwidth}{|
		>{\columncolor[HTML]{ECF4FF}}>{\centering\arraybackslash}X |>{\centering\arraybackslash}X|>{\centering\arraybackslash}X|>{\centering\arraybackslash}X|>{\centering\arraybackslash}X|}
\hline
& \cellcolor[HTML]{ECF4FF}\textbf{$a_{\rm c}$} & \cellcolor[HTML]{ECF4FF}\textbf{$\tau_{\rm c}$} & \cellcolor[HTML]{ECF4FF}\textbf{$\nu^{(i)}$} & \cellcolor[HTML]{ECF4FF}\textbf{$\nu^{(a)}$} \\ \hline
\textbf{\texttt{HH1}} & -0.68609                               & 0.19214                                   & 7.83301                                      & 8.51761                                      \\ \hline
\textbf{\texttt{HH2}} & -0.83431                               & 0.57621                                   & 8.58402                                      & 7.77241                                      \\ \hline
\textbf{\texttt{HH3}} & -1.33683                               & 0.95920                                   & 7.61733                                      & 8.75879                                      \\ \hline
\textbf{\texttt{HH4}} & -3.37162                               & 0.95457                                   & 7.23890                                      & 9.21666                                      \\ \hline
\end{tabularx}
\caption[Parameter values for the double Hopf bifurcations in the $(a,\tau)$ plane]{Parameter values for the double Hopf bifurcations in the $(a,\tau)$ plane, with $\gamma=0.641,\ \omega=2.6\pi$.}
\label{tab:Hopf_points}
\end{table*}


The normal form of the double Hopf bifurcation can be expressed in polar coordinates, with amplitudes $r_1$, $r_2$ \eqref{Hopf:doublenormalform_r} and phase angles $\varphi_1$, $\varphi_2$ \eqref{Hopf:doublenormalform_phi}; see also \cite{Molnar2017}. In \eqref{Hopf:doublenormalform_r}, $r_1$ corresponds to eigenvalue ${\rm i}\nu^{(i)}$, and $r_2$ corresponds to eigenvalue ${\rm i}\nu^{(a)}$. So, $r_1\neq0$, $r_2=0$ steady states of \eqref{Hopf:doublenormalform_r} correspond to in-phase limit cycles, and $r_1=0$, $r_2\neq0$ corresponds to anti-phase limit cycles.  

Table~\ref{tab:normal_form_coeffs} gives values for the normal form coefficients $a_{jk}$ and parameters $\rho_{jk}$, used in the expressions \eqref{Hopf:unfoldingparameters} for the unfolding parameters, for each of the four double Hopf points.
Note that the normal form coefficients can also be computed for a wide class of time delay systems numerically using the approach proposed in~\cite{Bosschaert2020} that was implemented as part of \textsc{DDE-Biftool}.
In the sequel, we explore the point \texttt{HH1} in detail.

\begin{table*}[t]
	\centering
	\def\arraystretch{1.5}
	\setlength{\arrayrulewidth}{1.5pt}
	\small
\begin{tabularx}{\textwidth}{|
		>{\columncolor[HTML]{ECF4FF}}c |>{\centering\arraybackslash}X|>{\centering\arraybackslash}X|>{\centering\arraybackslash}X|>{\centering\arraybackslash}X|>{\centering\arraybackslash}X|>{\centering\arraybackslash}X|>{\centering\arraybackslash}X|>{\centering\arraybackslash}X|}
\hline
& \cellcolor[HTML]{ECF4FF}\textbf{$a_{11}$} & \cellcolor[HTML]{ECF4FF}\textbf{$a_{12}$} & \cellcolor[HTML]{ECF4FF}\textbf{$a_{21}$} & \cellcolor[HTML]{ECF4FF}\textbf{$a_{22}$} & \cellcolor[HTML]{ECF4FF}\textbf{$\rho_{11}$} & \cellcolor[HTML]{ECF4FF}\textbf{$\rho_{12}$} & \cellcolor[HTML]{ECF4FF}\textbf{$\rho_{21}$} & \cellcolor[HTML]{ECF4FF}\textbf{$\rho_{22}$} \\ \hline
\textbf{\texttt{HH1}} & -1.45930                                  & -2.98167                                  & -3.17291                                  & -1.62071                                  & 0.41422                                      & -2.53782                                     & 0.44713                                      & 2.99651                                      \\ \hline
\textbf{\texttt{HH2}} & -1.57279                                  & -3.07152                                  & -2.88168                                  & -1.40684                                  & 0.25537                                      & 3.09705                                      & 0.24324                                      & -2.58847                                     \\ \hline
\textbf{\texttt{HH3}} & -1.05553                                  & -2.17940                                  & -2.28045                                  & -1.17705                                  & 0.04553                                      & -2.04607                                     & 0.03689                                      & 2.44805                                      \\ \hline
\textbf{\texttt{HH4}} & -0.61111                                  & -1.29100                                  & -1.31866                                  & -0.69641                                  & -0.00656                                     & -1.18880                                     & -0.01192                                     & 1.42407                                      \\ \hline
\end{tabularx}
	\caption[The values of the normal form coefficients and parameters $\rho_{jk}$ for the double Hopf bifurcations in the $(a,\tau)$ plane]{The values of the normal form coefficients and parameters $\rho_{jk}$, used in the expressions for the unfolding parameters, for each of the double Hopf bifurcation points in the $(a,\tau)$ plane. Fixed parameter values are $\alpha=12.457,\ \beta=0.007095,\ \gamma=0.641,\ \omega=2.6\pi,\ b=1$.}
	\label{tab:normal_form_coeffs}
\end{table*}


\subsection{Phase portraits}
\label{Hopf:subsec:PhasePortraits}
The normal form for the double Hopf bifurcation \eqref{Hopf:doublenormalform_r} gives rise to distinct structurally stable phase portraits in different regions of parameter space about the double Hopf points. To see this, we look at the possible steady states of \eqref{Hopf:doublenormalform_r},
\begin{equation}
\label{eq:amplitudes}
\begin{aligned}
\circled{1}:\ \ \ (r_1,r_2)&=(0,0),\\
\circled{2}:\ \ \ (r_1,r_2)&=\left(\sqrt{-\frac{b_1}{a_{11}}},0\right),\\
\circled{3}:\ \ \ (r_1,r_2)&=\left(0,\sqrt{-\frac{b_2}{a_{22}}}\right),\\
\circled{4}:\ \ \ (r_1,r_2)&=\left(\sqrt{\frac{a_{12}b_2-a_{22}b_1}{a_{11}a_{22}-a_{12}a_{21}}},
\sqrt{\frac{a_{21}b_1-a_{11}b_2}{a_{11}a_{22}-a_{12}a_{21}}}\right).
\end{aligned}
\end{equation}

According to Guckenheimer and Holmes \cite{Guckenheimer1983}, there is a partition of parameter space into regions with topologically different phase portraits in which different combinations of steady states exist simultaneously. The lines separating these regions are given by
\begin{equation}
\label{Hopf:lines1to4}
\begin{aligned}
\text{Line 1:}\ \ &b_1=0, \\
\text{Line 2:}\ \ &b_2=0, \\
\text{Line 3:}\ \ &a_{12}b_2-a_{22}b_1=0, \\
\text{Line 4:}\ \ &a_{21}b_1-a_{11}b_2=0.
\end{aligned}
\end{equation}

\begin{figure}
	\centering
	\includegraphics[width=\columnwidth]{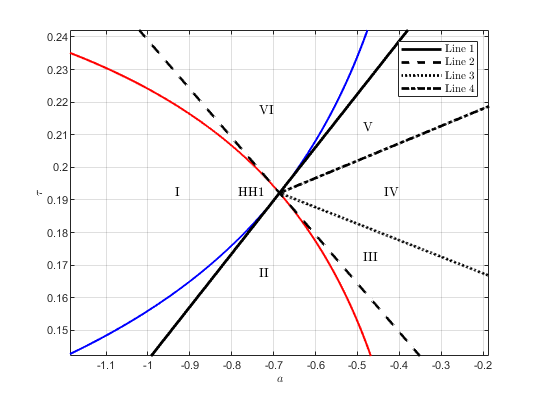}
	\caption[Regions with qualitatively different phase portraits about double Hopf point \texttt{HH1}]{Regions with qualitatively different phase portraits in the vicinity of double Hopf point \texttt{HH1}: $(a,\tau)=(-0.68609,0.19214)$. The blue curve corresponds to $\text{Re}(\lambda^{(i)})=0$ and the red curve to $\text{Re}(\lambda^{(a)})=0$. Lines 1-4 are defined by \eqref{Hopf:lines1to4}. Fixed parameter values are $\alpha=12.457,\ \beta=0.007095,\ \gamma=0.641,\ \omega=2.6\pi,\ b=1$.}
	\label{fig:Hopf:Hopf1_PhasePortraitRegions}
\end{figure}

Lines 1-4 divide parameter space into six regions, labelled I-VI in Figure \ref{fig:Hopf:Hopf1_PhasePortraitRegions}. There is excellent agreement near \texttt{HH1} between lines 1 and 2 in \eqref{Hopf:lines1to4} and the bifurcation curves, taken from Figure~\ref{fig:linear_stability_charts:Kuang_a_tau}.

The stability of the steady states \circled{1} -- \circled{4} is found by evaluating the Jacobian $J$ of \eqref{Hopf:doublenormalform_r}, given by,
\begin{align}
\begin{split}
J&=\begin{bmatrix}
J_{11} & J_{12} \\
J_{21} & J_{22}
\end{bmatrix},\\
J_{11}&=3a_{11}r_2^2+a_{12}r_2^2
+\rho_{11}\left(a-a_{\rm c}\right)+\rho_{12}\left(\tau-\tau_{\rm c}\right),\\
J_{12}&=J_{21}=2a_{21}r_1r_2,\\
J_{22}&=a_{21}r_1^2+3a_{22}r_2^2
+\rho_{21}\left(a-a_{\rm c}\right)+\rho_{22}\left(\tau-\tau_{\rm c}\right).
\end{split}
\end{align}

\begin{figure}
	\centering
	\begin{overpic}[width=.7\columnwidth]{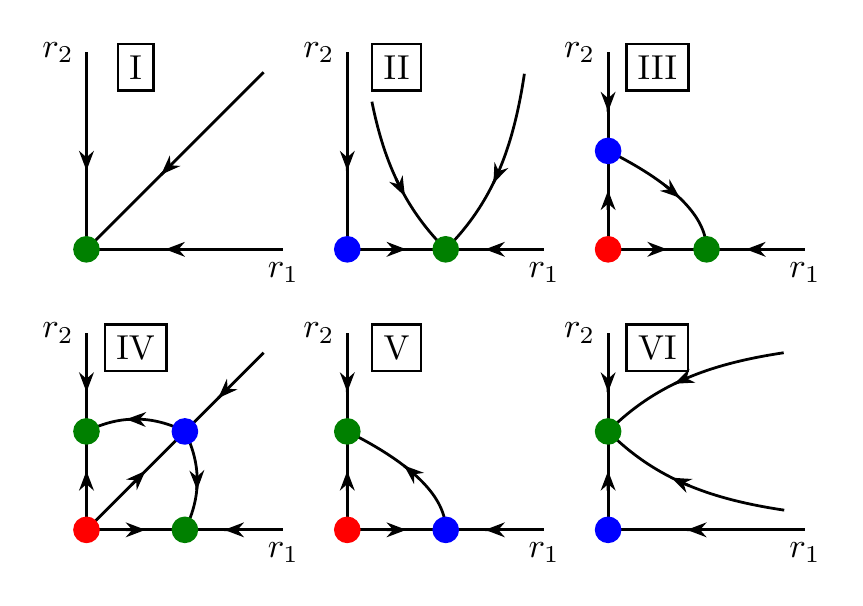}
       \put(8.5,37){{\fontsize{5}{6}$\circled{1}$}}
       
       \put(38.5,37){{\fontsize{5}{6}$\circled{1}$}}
       \put(50,37){{\fontsize{5}{6}$\circled{2}$}}
       
       \put(68.5,37){{\fontsize{5}{6}$\circled{1}$}}
       \put(80,37){{\fontsize{5}{6}$\circled{2}$}}
       \put(64.5,52){{\fontsize{5}{6}$\circled{3}$}}

       \put(8.5,4.5){{\fontsize{5}{6}$\circled{1}$}}
       \put(20,4.5){{\fontsize{5}{6}$\circled{2}$}}
       \put(4.5,19){{\fontsize{5}{6}$\circled{3}$}}
       \put(23.5,19){{\fontsize{5}{6}$\circled{4}$}}

       \put(38.5,4.5){{\fontsize{5}{6}$\circled{1}$}}
       \put(50,4.5){{\fontsize{5}{6}$\circled{2}$}}
       \put(34.5,19){{\fontsize{5}{6}$\circled{3}$}}

       \put(68.5,4.5){{\fontsize{5}{6}$\circled{1}$}}
       \put(64.5,19){{\fontsize{5}{6}$\circled{3}$}}  
    \end{overpic}
	\caption[Possible phase portrait topologies around double Hopf point \texttt{HH1}]{Possible phase portrait topologies around double Hopf point \texttt{HH1}. Sources, sinks and saddles are illustrated using red, green and blue dots respectively. Roman numerals correspond to the regions in Figure~\ref{fig:Hopf:Hopf1_PhasePortraitRegions}. Equilibria {\fontsize{5}{6}$\circled{1}-\circled{4}$} are given in \eqref{eq:amplitudes}.}
	\label{fig:Hopf:Hopf1_PhasePortraits}
\end{figure}

The analysis leads to the phase portraits in Figure~\ref{fig:Hopf:Hopf1_PhasePortraits}. 
Let us consider the qualitative changes in these phase portraits as we move counter-clockwise around the double Hopf point \texttt{HH1} in Figure \ref{fig:Hopf:Hopf1_PhasePortraitRegions}. In region I, only the zero equilibrium \circled{1} steady state exists. In agreement with analysis in Section~\ref{chap:linear_stability_charts}, this is stable. A Hopf bifurcation gives rise to the stable in-phase limit cycle \circled{2} in region II, and the zero equilibrium becomes a saddle. Moving into region III, an unstable anti-phase limit cycle \circled{3} is born and the zero equilibrium becomes a source. Region IV is a region of bistability of the in-phase and anti-phase limit cycles, together with an unstable quasi-periodic orbit \circled{4}. This quasi-periodic orbit then collides with the in-phase limit cycle \circled{2} to give the phase portrait in region V, where the anti-phase limit cycle \circled{3} is the only stable steady state. In region VI, the in-phase limit cycle no longer exists, leaving only the stable anti-phase limit cycle \circled{3} and the unstable equilibrium at the origin \circled{1}. This limit cycle disappears at a Hopf bifurcation as we move back into region I, and the equilibrium at the origin \circled{1} regains stability.   

\subsection{One-parameter bifurcation diagrams}
\begin{figure}[t]
\centering
\begin{subfigure}
\centering
\includegraphics[width=.9\columnwidth]{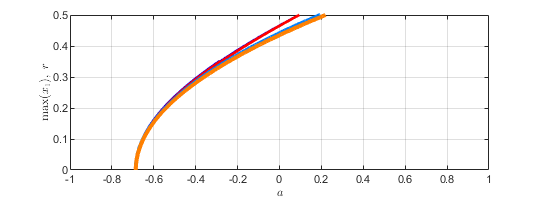}
\end{subfigure}
\begin{subfigure}
\centering
\includegraphics[width=.9\columnwidth]{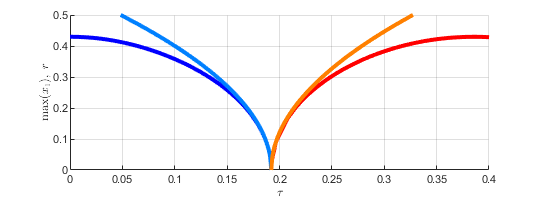}
\end{subfigure}
\begin{subfigure}
    \centering
    \includegraphics[width=.75\columnwidth]{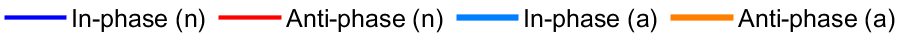}
\end{subfigure}
	\caption[One-parameter bifurcation diagrams in $a$ and $\tau$ for the periodic orbits, comparing analysis with numerics]{One-parameter bifurcation diagrams in $a$ (with $\tau=0.19214$ from Table~\ref{tab:Hopf_points}), and $\tau$ (with $a=-0.68609$) for the limit cycles born at \texttt{HH1}. Numerics (n) are compared with the analytic (a) solutions obtained using the normal form given in \eqref{Hopf:doublenormalform_r}. All solutions are stable at the bifurcation point. Fixed parameter values are $\alpha=12.457,\ \beta=0.007095,\ \gamma=0.641,\ \omega=2.6\pi,\ b=1$.}
	\label{fig:Hopf:oneparameterbifurcation}
 \end{figure}

Apart from obtaining different phase portraits around the double Hopf bifurcation, the normal form \eqref{Hopf:doublenormalform_r} allows us to plot one-parameter bifurcation diagrams to illustrate the unfolding of the in-phase and anti-phase solutions. These provide a useful means of comparison with numerical results. Figure~\ref{fig:Hopf:oneparameterbifurcation} illustrates the one-parameter bifurcation diagrams in $a$ (with $\tau=0.19214$ from Table~\ref{tab:Hopf_points}), and $\tau$ (with $a=-0.68609$) showing how the double Hopf bifurcation \texttt{HH1} unfolds. Numerical results obtained from \textsc{DDE-Biftool} are compared with the analytic results given by \eqref{eq:amplitudes}. 
The bifurcation diagrams in parameter $a$ show transitions between regions I and IV, and the birth of the two stable limit cycles. In contrast, the bifurcation diagrams in $\tau$ show transitions between regions II and VI, showing the destruction of one type of limit cycle, followed by the birth of the other. 

The centre manifold analysis is valid for small amplitudes near the double Hopf bifurcation point \texttt{HH1}.
Figure~\ref{fig:Hopf:oneparameterbifurcation} shows that there is strong agreement between the analytic and numeric solutions within the expected parameter range.  It only shows stable solutions. Other solutions will be discussed in the next section.

In Appendix~\ref{appC}, we give a brief discussion on the possibility, or otherwise, of internal resonances in this problem.

\section{Numerical results}
\label{chap:numerics}

\begin{figure}
	\centering
	\includegraphics[scale=1]{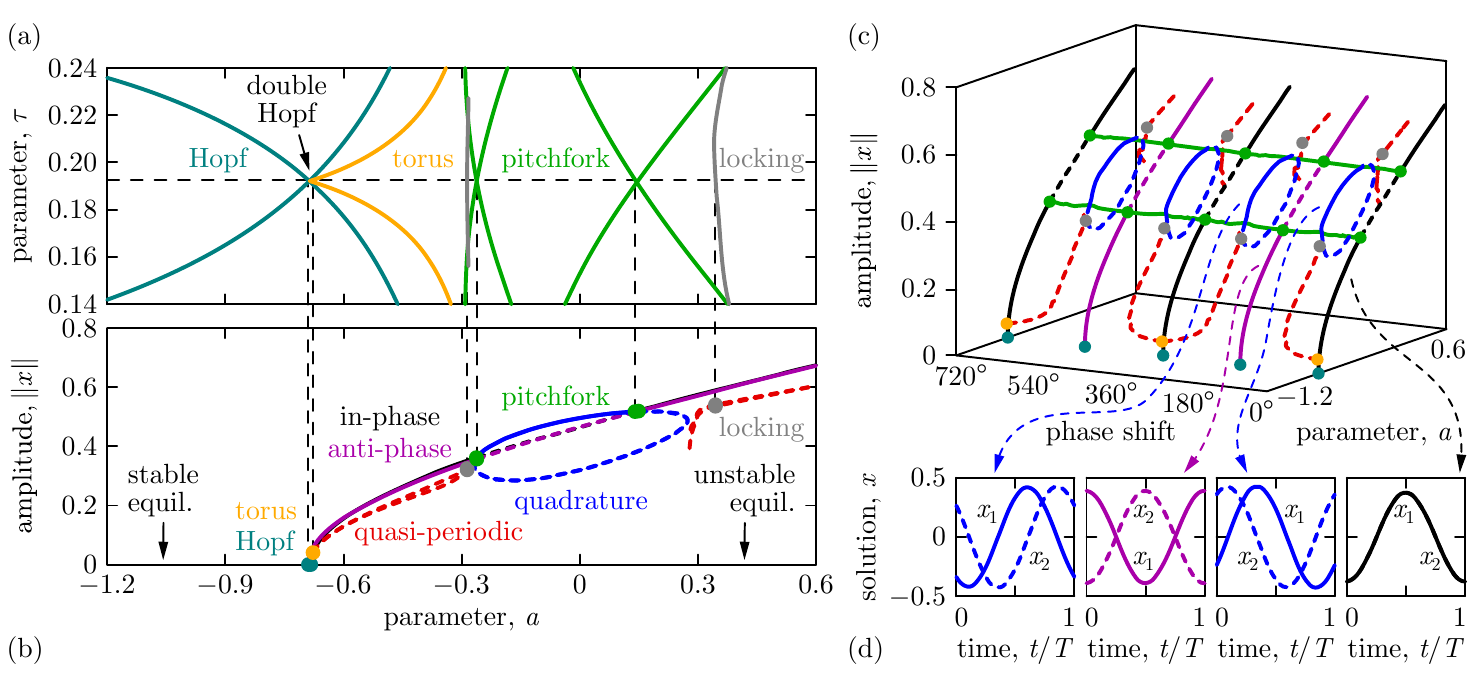}
	\caption[Numerical bifurcation diagrams of the delayed HKB system.]{Numerical bifurcation diagrams of the delayed HKB system \eqref{HKBeqns:full_delayed_HKB}.
	(a) Two-parameter diagram in the $(a,\tau)$ plane with branches of Hopf, torus, pitchfork bifurcation and 1-1 locking.
	(b) One-parameter diagram of limit cycles and quasi-periodic orbits against parameter $a$ for ${\tau=0.1926}$ (solid line: stable solution, dashed line: unstable solution).
	(c) One-parameter diagram also indicating the phase shift between $x_1$ and $x_2$.
	(d) In-phase, anti-phase, and phase quadrature periodic orbits at ${a=-0.2}$.
	Fixed parameter values are $\alpha=12.457,\ \beta=0.007095,\ \gamma=0.641,\ \omega=2.6\pi,\ b=1$.}
	\label{fig:continuation}
\end{figure}

To further analyse the global behaviour of the nonlinear delayed HKB model \eqref{HKBeqns:full_delayed_HKB}, we conducted extensive numerical bifurcation calculations\footnote{Code used to obtain these results is available at \url{https://github.com/DomboZoli/Quasi-periodic-package}.}.
The results are summarised in Figure \ref{fig:continuation}.
We computed branches of equilibria and limit cycles using \textsc{DDE-Biftool} \cite{DDEBifTool2002} and branches of quasi-periodic orbits using the Matlab tool introduced in \cite{Dombovari2015} and further developed in \cite{Molnar2017}.

We analysed the stability of the trivial equilibrium, detected Hopf bifurcations, and continued the branches of Hopf bifurcation in two parameters, $a$ and $\tau$; see the teal branches\footnote{These are the red and blue branches in Figure~\ref{fig:Hopf:Hopf1_PhasePortraitRegions}.} in Figure \ref{fig:continuation}(a).
The two branches give rise to the in-phase and anti-phase limit cycles, and their the intersection is the double-Hopf bifurcation point \texttt{HH1} at $(a,\tau)=(a_c,\tau_c)=(-0.68609,0.19214)$.
Then, we selected a delay value ${\tau=0.1926}$ close to \texttt{HH1} and continued the in-phase and anti-phase limit cycles by varying parameter $a$; these are the black and purple curves, respectively, in Figure \ref{fig:continuation}(b), which lie almost on top of one another and above the red quasi-periodic curve\footnote{The amplitude measure on the vertical axis is the ``root-mean-square value'' over the period $T$: ${\| x \| = \sqrt{\frac{1}{T} \int_{0}^{T} (x_1^2(t) + x_2^2(t)) {\rm d} t}}$.}.

Then, we considered the stability of the in-phase and anti-phase limit cycles. We detected torus bifurcations (orange points) and pitchfork bifurcations (green points), associated with a pair complex and one real characteristic multipliers located on the unit circle of the complex plane, respectively.
These bifurcations were continued in two parameters as shown by the orange and green branches in Figure \ref{fig:continuation}(a).

The pitchfork bifurcations give rise to additional limit cycles.
These solutions are plotted in blue in Figure \ref{fig:continuation}(c) as a function of the phase shift\footnote{Based on the location of the maximum points of $x_1(t)$ and $x_2(t)$.} between the states $x_1$ and $x_2$.
The phase shift of the limit cycles changes continuously and sweeps across the entire $[0^\circ,360^\circ]$ domain.
The diagram repeats every $360^\circ$ along the phase shift axis.

When the phase shift of these additional limit cycles reaches approximately $90^\circ$ or $270^\circ$, the branch splits into circular branches of limit cycles, where the associated phase shift was observed to be approximately constant $90^\circ$ or $270^\circ$, respectively.
We refer to these as \textit{limit cycles in phase quadrature}.
The circular branches of these solutions can be seen in Figure \ref{fig:continuation}(b). The in-phase, anti-phase and phase quadrature limit cycles themselves are depicted in Figure \ref{fig:continuation}(d) for ${a=-0.2}$.
For these phase quadrature limit cycles, we further observed that the time period is about 4 times the delay, and therefore a special type of symmetry with $x_1(t) \approx x_2(t-\tau)$ or $x_2(t) \approx x_1(t-\tau)$ holds, respectively.

The torus bifurcations (orange points in Figure \ref{fig:continuation}) give rise to two branches of quasi-periodic orbits (red dashed curves)\footnote{It is not possible to use \textsc{DDE-Biftool} to automatically continue the quasi-periodic orbit from the torus bifurcation that occurs when a limit cycle changes stability. Domb\'ov\'ari and St\'ep\'an \cite{Dombovari2015} point out that saddle-like invariant sets are especially difficult to find for    DDEs because the standard trick of tracking solutions along reversed time cannot be used. Their algorithm is not straightforward to implement, and success is dependent on the accuracy of the initial solution profile estimates and on stability properties of the solution in question.}.
The two branches are identical except that $x_1$ and $x_2$ are interchanged\footnote{Note that system~(\ref{HKBeqns:full_delayed_HKB}) is symmetric in $x_1$ and $x_2$, whereas the quasi-periodic orbits branch out from the in-phase solution where ${x_1(t) \equiv x_2(t)}$.}.

The quasi-periodic orbits are associated with two angular frequencies, $\upsilon_1$ and $\upsilon_2$, that are close to the frequencies related to the double Hopf point \texttt{HH1}.
Accordingly, quasi-periodic orbits are described as a surface (a torus) parameterised by two dimensionless time variables $\theta_1 \in [0,2\pi]$ and $\theta_2 \in [0,2\pi]$ associated with $\upsilon_1$ and $\upsilon_2$.
The quasi-periodic branches are indicated with their phase shift\footnote{The amplitude and phase shift were determined based on taking a section $\theta_1={\rm mod}(\upsilon_1 \theta/\upsilon_2,2\pi)$, $\theta_2={\rm mod}(\upsilon_2 \theta/\upsilon_1,2\pi)$ of the torus with ${\theta \in [0,2\pi]}$, and then using the same amplitude and phase measures for $x_1(t)$ and $x_2(t)$ as for limit cycles with period ${T=2\pi}$.} in Figure \ref{fig:continuation}(c).

Finally, as the quasi-periodic branch is continued, multiple 1-1 locking states, coloured gray in  Figure \ref{fig:continuation}, were detected where ${\upsilon_1/\upsilon_2=1}$.
At the locking points, the quasi-periodic orbit degrades into the limit cycle in phase quadrature.
Near locking, the corresponding torus becomes challenging to compute numerically; see the loss of accuracy along the red curve on the right of Figure \ref{fig:continuation}(b).

An approximation of the locking point was continued in the $(a,\tau)$ plane by establishing a numerical condition on higher harmonics along both dimensionless time coordinates. These 
harmonics were determined by discrete Fourier transform (DFT), in which the corresponding derivatives with respect to present states, retarded states and parameters were determined analytically to give a well conditioned two parameter continuation scheme.
This continuation scheme is detailed in Appendix~\ref{appD}.
The numerical locking points lie closer to the branch of limit cycles in phase quadrature (blue curve) when $a=-0.29$ than for when $a=0.34$.

The results in this section provide further information about the global dynamic behaviour of the HKB system \eqref{HKBeqns:full_delayed_HKB}, that could not be obtained by the analytic methods presented in previous sections.

\section{Conclusions}
In this paper, we discussed the effects of delay on the Haken-Kelso-Bunz (HKB) model \cite{Haken1985} of bimanual human motor coordination.
We investigated the stability of the trivial solutions in the corresponding linear system \eqref{HKB:normal_forms_together}, which can be written as a delayed oscillator where the feedback changes the damping. We discovered
Hopf and double Hopf bifurcations in this linear delayed HKB model. We analysed the double Hopf bifurcations in the full HKB system \eqref{HKBeqns:full_delayed_HKB} by means of centre manifold reduction \cite{Haken1985} to calculate the stability of both in-phase and anti-phase limit cycles and quasi-periodic orbits. We verified our results using numerical continuation. In addition, we discovered limit cycles in phase quadrature and 1-1 locking of quasi-periodic orbits. We have shown that in-phase and the anti-phase limit cycles can be replaced by a phase-lagged solution via a pitchfork bifurcation of periodic orbits. This phase lagged solution has transitioned from the quasi-periodic branch which emerges from the double Hopf bifurcation point. This suggests that double frequency transient transitioning is needed to change from in-phase to anti-phase limit cycles, while the phase-lagged solution is reached by increasing the linear coupling coefficient $a$.
The results provide valuable insights into the nonlinear dynamic behaviour of this model,
which may help determine the relevance of the delayed HKB system to the application of the mirror game for early diagnosis of disorders such as schizophrenia~\cite{Varlet2012}. Furthermore, we suggest that the methods we have presented may be valuable to assess the corresponding suitability of other delayed models of human motor coordination, such as the variation of the delayed HKB system presented by S\l{}owi\'{n}ski et al.~\cite{Slowinski2020}, which incorporates a neurologically motivated coupling term.


\appendix
\section{Stability charts}
\label{appA}
Equation \eqref{HKB:normal_forms_together_control_form} has the trivial solution $\eta(t)=0$. We want to find conditions under which this solution is stable. If this solution becomes unstable via a Hopf bifurcation, we would expect to find finite amplitude limit cycles, corresponding to observable oscillations in the full HKB system \eqref{HKBeqns:full_delayed_HKB}. The characteristic equation of \eqref{HKB:normal_forms_together_control_form} is given by
\begin{equation}
\label{linear_stability_charts:HKBchareqn}
\lambda^2-(\gamma+a)\lambda+\omega^2=\mp a\lambda {\rm e}^{-\lambda\tau}.
\end{equation}
The solution $\eta(t)=0$ of \eqref{HKB:normal_forms_together_control_form} is stable when ${\text{Re}(\lambda)<0}$. 

When $\tau\neq 0$, \eqref{linear_stability_charts:HKBchareqn} is an exponential polynomial in $\lambda$, which has an infinite number of roots (either real or complex conjugate). If {\it any} of these roots have positive real part, then the steady state $\eta(t)=0$ is unstable. 

From Kuang \cite[Theorem~1.4, p.~66 \& Section 3.3]{Kuang1993}, we know that stability changes (or boundaries) occur when the root with the largest real part is purely imaginary. The locations in parameter space where such roots exist can be found by substituting $\lambda={\rm i}\nu$, where $\nu\ge0$ is real, into \eqref{linear_stability_charts:HKBchareqn} and equating real and imaginary parts, to give
\begin{equation}
\label{linear_stability_charts:HKBlin_charri}
\begin{aligned}
\omega^2-\nu^2\pm a\nu\sin(\nu\tau)&=0,\\
-(\gamma+a)\nu\pm a\nu\cos(\nu\tau)&=0.
\end{aligned}
\end{equation}    
Note that $\nu=0$ is a solution of \eqref{linear_stability_charts:HKBlin_charri} when $\omega=0$. Hence we take $\nu>0$ in the sequel. From \eqref{linear_stability_charts:HKBlin_charri} we obtain a quartic expression in $\nu$,
\begin{equation}
\label{linear_stability_charts:Kuang_nupoly}
\nu^4+\left(\left(\gamma+a\right)^2-a^2-2\omega^2\right)\nu^2+\omega^4=0,
\end{equation}
which has roots
\begin{equation}
\label{linear_stability_charts:nu_pm}
\begin{aligned}
\nu_{\pm}^2&=\frac{1}{2}\Big(a^2+2\omega^2-(\gamma+a)^2
\pm\sqrt{(a^2+2\omega^2-(\gamma+a)^2)^2-4\omega^4}\Big).
\end{aligned}
\end{equation}
Since $\nu_{\pm}$ must be real, we have $\lambda_{\pm}={\rm i}\nu_{\pm}$, $\nu_{+}>\nu_{-}>0$ provided that
\begin{equation}
\label{linear_stability_charts:conditionImRootmapped}
\begin{aligned}
(\text{a})\hspace{1cm} &a^2+2\omega^2-(\gamma+a)^2>0,\\
(\text{b})\hspace{1cm} &(a^2+2\omega^2-(\gamma+a)^2)^2>4\omega^4,
\end{aligned}
\end{equation}
with no solutions otherwise.

The real parts of the rightmost eigenvalues either become positive (instability) or negative (stability) when parameters such as $\tau$ change. But we do not know in which direction the eigenvalues move. To determine this direction, the sign of the derivative of $\text{Re}\left(\lambda\left(\tau\right)\right)$ with respect to $\tau$ needs to be found at the point where $\lambda\left(\tau\right)$ is purely imaginary. So we calculate the sign of $\rho_{\tau}$, given by
\begin{equation}
    \label{rho_definition}
    \begin{aligned}
    \rho_{\tau}:=\text{Re}\left(\frac{{\rm d}\lambda}{{\rm d}\tau}\biggr\rvert_{\lambda ={\rm i}\nu}\right).
    \end{aligned}
\end{equation}

It turns out to be more convenient to calculate the inverse of $\rho_{\tau}$, as only its sign matters. By differentiating the characteristic equation \eqref{linear_stability_charts:HKBchareqn}, we obtain
\begin{equation}
\label{linear_stability_charts:Kuang_dldt}
\begin{aligned}
&  \left(\frac{{\rm d}\lambda}{{\rm d}\tau}\right)^{-1}=\frac{\pm\left(2\lambda-\left(\gamma+a\right)\right){\rm e}^{\lambda\tau}+ a}{ a \lambda^2}-\frac{\tau}{\lambda}.
\end{aligned}
\end{equation}
Then substituting the expression for ${\rm e}^{\lambda\tau}$ obtained from \eqref{linear_stability_charts:HKBchareqn} into \eqref{linear_stability_charts:Kuang_dldt}, and using \eqref{linear_stability_charts:nu_pm}, we have
\begin{equation}
\begin{aligned}
&\text{sign}\rho_{\tau}=\text{sign}\left[\text{Re}\left(\left(\frac{{\rm d}\lambda}{{\rm d}\tau}\right)^{-1}\biggr\rvert_{\lambda ={\rm i}\nu_{\pm}}\right)\right]
=\text{sign}\left(\pm\sqrt{(a^2+2\omega^2-(\gamma+a)^2)^2-4\omega^4}\right).
\end{aligned}
\end{equation}
Thus $\rho_{\tau}>0$ (eigenvalues crossing the imaginary axis from left to right with increasing $\tau$) occurs for $\tau$ corresponding to $\nu_{+}$ and $\rho_{\tau}<0$ (eigenvalues crossing the imaginary axis from right to left with increasing $\tau$) occurs for $\tau$ corresponding to $\nu_{-}$. We find these values of $\tau$ from \eqref{linear_stability_charts:HKBlin_charri}, by setting 
\begin{equation}
\label{linear_stability_charts:Kuang_tau}
\begin{aligned}
\tau_{n,1}&=\frac{\theta_1}{\nu_{+}}+\frac{2n\pi}{\nu_{+}}\\
\tau_{n,2}&=\frac{\theta_2}{\nu_{-}}+\frac{2n\pi}{\nu_{-}}
\end{aligned}
\end{equation} 
for $n \in \mathbb{Z}$ where $\theta_{1,2} \in [0, 2\pi)$ are given by
\begin{equation}
\label{linear_stability_charts:Kuang_sincos}
\begin{split}
\cos(\theta_{1})&=\pm \frac{\gamma+a}{a},\\
\sin(\theta_{1})&=\pm \frac{\nu_+^2-\omega^2}{a\nu_+},
\end{split}
\quad
\begin{split}
\cos(\theta_{2})&=\pm \frac{\gamma+a}{a},\\
\sin(\theta_{2})&=\pm \frac{\nu_-^2-\omega^2}{a\nu_-}.
\end{split}
\end{equation}

Kuang's theorem \cite[Theorem~1.4, p.~66]{Kuang1993} allows us to consider the eigenvalues $\lambda$ as a continuous function of $\tau$. Then the stability of the solution $\eta(t)=0$ of \eqref{HKB:normal_forms_together_control_form} for $\tau>0$ can be found by looking at the stability of the system at $\tau=0$. 

If $\eta(t)=0$ is {\it stable} when $\tau=0$, then $\tau_{0,1}<\tau_{0,2}$ because the multiplicity of roots with positive real part cannot become negative. Additionally, we have that
\begin{equation}
\label{linear_stability_charts:Kuang_taudiff}
\tau_{n+1,1}-\tau_{n,1}=\frac{2\pi}{\nu_{+}}<\frac{2\pi}{\nu_{-}}=\tau_{n+1,2}-\tau_{n,2}.
\end{equation} 
This means that there can only be a finite number of switches between stability and instability. Specifically, there are $k$ switches from stability to instability to stability when 
\begin{equation}
\label{linear_stability_charts:Kuang_taustab}
\begin{aligned}
&\tau_{0,1}<\tau_{0,2}<\tau_{1,1}<\cdots
<\tau_{k-1,1}<\tau_{k-1,2}<\tau_{k,1}<\tau_{k+1,1}<\tau_{k,2}<\cdots.
\end{aligned}
\end{equation}

If $\eta(t)=0$ is {\it unstable} when $\tau=0$, then it is either unstable for $\tau>0$, or a finite number of stability switches occur; $k$ switches from instability to stability to instability \textit{may} occur when
\begin{equation}
\label{linear_stability_charts:Kuang_tauunstab}
\begin{aligned}
&\tau_{0,2}<\tau_{0,1}<\tau_{1,2}<\cdots
<\tau_{k-1,2}<\tau_{k-1,1}<\tau_{k,1}<\tau_{k,2}<\cdots\ .
\end{aligned}
\end{equation} 
%
It is straightforward to show that this condition is satisfied for \eqref{HKB:normal_forms_together_control_form}. 
This means that stability changes occur according to \eqref{linear_stability_charts:Kuang_tauunstab} when the system is unstable at $\tau=0$.

Stability charts in parameter spaces $(a, \tau)$ and $(\gamma,\tau)$ (Figures~\ref{fig:linear_stability_charts:Kuang_a_tau} and \ref{fig:linear_stability_charts:Kuang_gamma_tau}) are obtained by treating \eqref{linear_stability_charts:Kuang_tau} and \eqref{linear_stability_charts:Kuang_sincos} as functions of $a$ and $\gamma$, respectively, and carefully checking the conditions in \eqref{linear_stability_charts:Kuang_taustab} and \eqref{linear_stability_charts:Kuang_tauunstab} to see if the resulting curves correspond to stability changes. 

The above analysis has to be extended when we examine stability in $(\gamma,a)$ parameter space for fixed $\tau$. We adapt Kuang's \cite{Kuang1993} analysis to derive a parameterisation of the $\text{Re}(\lambda)=0$ curves directly. Rearranging \eqref{linear_stability_charts:HKBlin_charri}, we get expressions for $a$ and $\gamma$ in terms of $\nu$ as
\begin{equation}
\label{linear_stability_charts:HKBlin_gamma_a_parametrisation}
\begin{aligned}
a&=\pm \frac{\nu^2-\omega^2}{\nu\sin(\nu\tau)},\\
\gamma&= \frac{\nu^2-\omega^2}{\nu\sin(\nu\tau)} \left(\cos(\nu\tau)\mp 1\right),
\end{aligned}
\end{equation}
for $\nu\tau\neq n\pi$, $n\in\mathbb{Z}$. 

When $\nu\tau = n\pi$, $n\in\mathbb{Z}$, we have $\nu=\omega$ from the first equation of \eqref{linear_stability_charts:HKBlin_charri}. The second equation then gives \begin{equation}
\label{linear_stability_charts:HKBlin_Dcurve_special_case}
\begin{aligned}
& \gamma=a(\pm(-1)^n-1) \ \ \text{for $\omega\neq0$}.
\end{aligned}
\end{equation}
We now return to \eqref{linear_stability_charts:HKBchareqn} and look at the sign of $\text{Re}\left(\left(\frac{{\rm d}\lambda}{d\gamma}\right)^{-1}\right)$ evaluated on the curves given by \eqref{linear_stability_charts:HKBlin_gamma_a_parametrisation} and \eqref{linear_stability_charts:HKBlin_Dcurve_special_case}.  
We find that
\begin{equation}
\label{linear_stability_charts:HKBlin_sign_eqn}
\begin{aligned}
\text{Re}\left(\left(\frac{{\rm d}\lambda}{d\gamma}\right)^{-1}\biggr |_{\lambda={\rm i}\nu}\right)&=1+\frac{\omega^2}{\nu^2}-(\gamma+a)\tau =: \xi.
\end{aligned}
\end{equation}
In addition it is useful to look at $\frac{{\rm d}a}{{\rm d}\gamma}$ to see how the change of sign of \eqref{linear_stability_charts:HKBlin_sign_eqn} relates to the curves given by \eqref{linear_stability_charts:HKBlin_gamma_a_parametrisation}.
We find
\begin{equation}
\begin{aligned}
\frac{{\rm d}a}{{\rm d}\gamma}
&=\frac{\xi}{\frac{\gamma}{a}\xi-a\tau\sin^2(\nu\tau)},
\end{aligned}
\end{equation}
where $\xi$ is defined in \eqref{linear_stability_charts:HKBlin_sign_eqn}. Hence $\frac{{\rm d}a}{{\rm d}\gamma}=0$ if and only if $\xi=0$, for $a\neq0$. Therefore, the points in parameter space where the eigenvalues change direction correspond to the turning points of the $\text{Re}(\lambda)=0$ curves in $(\gamma,a)$ space. 
\section{Centre manifold reduction}
\label{appB}
First we will outline the method for calculating the normal form of a generic Hopf bifurcation of a DDE system using centre manifold theory and then apply this to find the normal form of the double Hopf bifurcations of the delayed HKB system. 
The approach given here is outlined by several authors \cite{Campbell2009, Kalmar2001, Molnar2017}). The detailed theory is discussed by Hale and Verduyn Lunel \cite{Hale1993}. 

This analysis applies to retarded delay differential equations with constant delay $\tau>0$. Consider a general delay differential equation of this type,
\begin{equation}
\label{Hopf:generalDDE}
\dot{\mathbf{x}}(t)=\mathbf{g}(\mathbf{x}(t),\mathbf{x}(t-\tau);\mu),
\end{equation}
where $\mathbf{x}\in\mathbb{R}^n$, $\mathbf{g}:\mathbb{R}^n\times\mathbb{R}^n\times\mathbb{R}^k\rightarrow\mathbb{R}^n$, $n,k \in \mathbb{Z}^+$, and $\mu\in\mathbb{R}^k$ and $\tau>0$ are parameters in the model. We assume that $\mathbf{g}$ is sufficiently smooth for the required computations and that the equation admits an equilibrium solution $\mathbf{x}_*$ which is independent of $\tau$. By shifting the equilibrium to zero and separating linear and nonlinear terms, \eqref{Hopf:generalDDE} can be written in the form
\begin{equation}
\label{Hopf:generalDDEtransformed}
\dot{\mathbf{x}}(t)=A_0(\mu)\mathbf{x}(t)+A_1(\mu)\mathbf{x}(t-\tau)\\+\mathbf{f}(\mathbf{x}(t),\mathbf{x}(t-\tau);\mu),
\end{equation}
where $A_j(\mu)=D_{j+1}\mathbf{g}(\mathbf{x}_*,\mathbf{x}_*;\mu)$ is the Jacobian of $\mathbf{g}$ with respect to its $(j+1)^{\rm th}$ argument, and 
\begin{equation}
\mathbf{f}(\mathbf{x}(t),\mathbf{x}(t-\tau);\mu)=\mathbf{g}(\mathbf{x}(t),\mathbf{x}(t-\tau);\mu)\\-A_0(\mu)\mathbf{x}(t)-A_1(\mu)\mathbf{x}(t-\tau).
\end{equation}

The characteristic equation of \eqref{Hopf:generalDDE} is then given by
\begin{equation}
\label{Hopf:chareqn}
\det(\Delta(\lambda;\mu))\\=\det(\lambda I_{n\times n} -A_0(\mu)-A_1(\mu){\rm e}^{-\lambda\tau})=0,
\end{equation}
where $I_{n\times n}$ is the $n\times n$ identity matrix.

The following analysis applies to critical parameter values $\mu=\mu_{\rm c}$ where the characteristic equation \eqref{Hopf:chareqn} has $m>0$ roots with zero real part, and the rest of the eigenvalues have negative real parts. We assume that the eigenvalues with zero real part have multiplicity one, which covers single and double Hopf bifurcations. 

To make progress with centre manifold construction, the operator differential equation representation of the DDE is required. Writing \eqref{Hopf:generalDDEtransformed} as an evolution equation on the Banach space $\mathcal{B}$ of continuously differentiable functions from $[-\tau,0]$ to $\mathbb{R}^n$ gives
\begin{equation}
\label{Hopf:operatorDE}
\dot{\mathbf{x}}_t=\mathcal{A}\mathbf{x}_t+\mathcal{F}(\mathbf{x}_t),
\end{equation}
where $\mathbf{x}_t\in\mathcal{B}$ is defined by
\begin{equation}
\label{Hopf:DDEflow}
\mathbf{x}_t(\theta)=\mathbf{x}(t+\theta),\quad \theta\in[-\tau,0],
\end{equation}
the linear operator $\mathcal{A}$ is defined by
\begin{align}
\label{Hopf:operatorA}
\mathcal{A}\phi(\theta)=
\begin{cases}
\frac{{\rm d}}{{\rm d}\theta}\phi(\theta),&\theta\in[-\tau,0),\\
A_0(\mu_{\rm c})\phi(0) +A_1(\mu_{\rm c})\phi(-\tau),
&\theta=0,
\end{cases}
\end{align}
and the nonlinear operator is
\begin{align}
\label{Hopf:nonlinearF}
\mathcal{F}(\phi)(\theta)=
\begin{cases}
0,& \theta\in[-\tau,0),\\
\mathbf{f}(\phi(0),\phi(-\tau);\mu_{\rm c}),&\theta=0.
\end{cases}
\end{align}

For the nonlinear calculations, it will be useful to define the operators
\begin{align}
\label{L,F}
\begin{split}
L(\phi)&=A_0(\mu_{\rm c})\phi(0)+A_1(\mu_{\rm c})\phi(-\tau),\\
\mathbf{F}(\phi)&=\mathbf{f}(\phi(0),\phi(-\tau);\mu_{\rm c}).
\end{split}
\end{align}

The following calculations will also require the dual space $\mathcal{B}^*$ of continuously differentiable functions on $[0,\tau]$ to $\mathbb{R}^{n*}$ (the $n$-dimensional row vectors), an adjoint operator
\begin{align}
\label{Hopf:adjointA}
\mathcal{A}^*\psi(\xi)=
\begin{cases}
-\frac{d}{d\xi}\psi(\xi),&\xi\in(0,\tau],\\
\psi(0)A_0(\mu_{\rm c})+\psi(\tau)A_1(\mu_{\rm c}),
&\xi=0,
\end{cases}
\end{align}
where we have assumed $A_0(\mu_{\rm c})$ and $A_1(\mu_{\rm c})$ are real, and the bilinear form $(\ ,\ ):\mathcal{B}^*\times\mathcal{B}\rightarrow\mathbb{R}$ given by
\begin{equation}
\label{Hopf:bilinearform}
(\psi,\phi)=\psi(0)\phi(0)+\int_{-\tau}^{0}\psi(\sigma+\tau)A_1(\mu_{\rm c})\phi(\sigma)d\sigma.
\end{equation}

The adjoint operator and bilinear form allow a projection of the solution to the DDE at the critical parameter values onto the centre manifold to be constructed. As with the ODE case, a first order approximation is constructed by considering the linear problem. Here the solution space can be decomposed as $\mathcal{B}=\mathcal{C}\bigoplus\mathcal S$ where $\mathcal{C}$ is an $m$-dimensional solution space spanned by the solutions corresponding to the eigenvalues with zero real part, $\mathcal{S}$ is infinite dimensional, and both $\mathcal{C}$ and $\mathcal{S}$ are invariant under the flow of the linear system. These are analogous to the centre and stable eigenspaces for ODEs. 

Let $\{\phi_1,\phi_2,\ldots,\phi_m\}$ be the basis for $\mathcal{C}$, with corresponding eigenvalues $\{\lambda_1,\lambda_2,\ldots,\lambda_m\}$. \linebreak Note that the eigenvalues of $\mathcal{A}$ are the same as the roots of the characteristic equation given by \eqref{Hopf:chareqn}. It was shown in Appendix~\ref{appA} that ${\lambda_k={\rm i}\nu_k}$ with $\nu_k\neq0$, therefore attention will be restricted to the case where all of the eigenvalues take this form. If $\lambda_k={\rm i}\nu_k$ is a root of \eqref{Hopf:chareqn}, then so is $-{\rm i}\nu_k$. The eigenvalues are labelled so that $\lambda_{k+1}=-{\rm i}\nu_k$, where $\nu_k>0$ and $k$ is odd. 

We find the basis for the centre eigenspace as follows. Consider a complex eigenfunction $\Phi_k\in\mathcal{B}$ corresponding to eigenvalue ${\rm i}\nu_k$ which satisfies $\forall \theta \in [-\tau,0]$ that
\begin{equation}
\mathcal{A}\Phi_k(\theta)={\rm i}\nu_k\Phi_k(\theta).
\end{equation}
Separating into real and imaginary parts gives
\begin{align}
\begin{split}
\mathcal{A}\phi_k(\theta)&=-\nu_k\phi_{k+1}(\theta),\\
\mathcal{A}\phi_{k+1}(\theta)&=\nu_k\phi_{k}(\theta),
\end{split}
\end{align}
where $\Phi_k(\theta)=\phi_k(\theta)+{\rm i}\phi_{k+1}(\theta)$.
For convenience, the basis shall be written as an $n\times m$ matrix, defined as
\begin{equation}
\label{Hopf:basis}
\bm{\Phi}(\theta)=
\begin{bmatrix}
\phi_1(\theta) & \phi_2(\theta) & \cdots & \phi_m(\theta)
\end{bmatrix}.
\end{equation}

Using the definition of $\mathcal A$, given by \eqref{Hopf:operatorA}, it can be shown that
\begin{equation}
\label{Hopf:PhiODE}
\bm{\Phi}'(\theta)=\bm{\Phi}(\theta)B,
\end{equation}
where $B$ is a block diagonal $m \times m$ matrix with blocks
\begin{equation}
\label{Hopf:Bblock}
B_k =
\begin{bmatrix}
0 & \nu_k \\
-\nu_k & 0
\end{bmatrix}
\end{equation}
for every pair of complex conjugate eigenvalues $\pm {\rm i}\nu_k$. It also follows from the definition of $\mathcal{A}$ that
\begin{equation}
\label{Hopf:boundaryeqn}
A_0(\mu_{\rm c})\bm{\Phi}(0)+A_1(\mu_{\rm c})\bm{\Phi}(-\tau)=\bm{\Phi}(0)B.
\end{equation} 
Solving \eqref{Hopf:PhiODE} with boundary condition \eqref{Hopf:boundaryeqn} gives
\begin{align}
\label{Hopf:solvePhiODE}
\begin{split}
\phi_k(\theta)&=\text{Re}({\rm e}^{{\rm i}\nu_k\theta}\mathbf{v}_k),\\
\phi_{k+1}(\theta)&=\text{Im}({\rm e}^{{\rm i}\nu_k\theta}\mathbf{v}_k),
\end{split}
\end{align}
where $\mathbf{v}_k$ satisfies $\Delta({\rm i}\nu_k;\mu_{\rm c})\mathbf{v}_k=0$, and $\Delta(;)$ is given in \eqref{Hopf:chareqn}.

The basis
\begin{equation}
\label{Hopf:adjointbasis}
\bm{\Psi}(\xi)=\begin{bmatrix}
\psi_1(\xi)\\
\vdots \\\psi_m(\xi)
\end{bmatrix}
\end{equation}
for the adjoint can be found in a similar way by deriving and solving 
\begin{align}
\begin{split}
\label{Hopf:adjointeqns}
& \bm{\Psi}'(\xi)=B\bm{\Psi}(\xi)\\
& \bm{\Psi}(0)A_0(\mu_{\rm c})+\bm{\Psi}(\tau)A_1(\mu_{\rm c})=-B\bm{\Psi}(0).
\end{split}
\end{align}

Note that the construction of $\mathcal{A}^*$ ensures that the eigenvalues of $\mathcal{A}^*$ are the same as the eigenvalues of $\mathcal{A}$, which is made explicit by considering eigenfunctions of the form $\psi(\xi)=\mathbf{w}{\rm e}^{-\lambda\xi}$, $\mathbf{w}\in\mathbb{R}^{n*}$. Here we label the eigenfunction in $\mathcal{B}$ with eigenvalue ${\rm i}\nu$ in the same way as the eigenfunction in $\mathcal{B}^{*}$ with the same eigenvalue.

Solving \eqref{Hopf:adjointeqns} yields
\begin{align}
\label{Hopf:solvePsiODE}
\begin{split}
\psi_k(\xi)&=\text{Re}(\mathbf{w}_k{\rm e}^{-{\rm i}\nu_k\xi}),\\
\psi_{k+1}(\xi)&=\text{Im}(\mathbf{w}_k{\rm e}^{-{\rm i}\nu_k\xi}),
\end{split}
\end{align}
where $\mathbf{w}_k\Delta({\rm i}\nu_k;\mu_{\rm c})=0$. Using the remaining degrees of freedom, the bases can be chosen such that $(\bm{\Psi},\bm{\Phi})=I_{m\times m}$, where $(\bm{\Psi},\bm{\Phi})$ is the matrix with $i,j$ elements $(\psi_i,\phi_j)$. Note that $\bm{\Psi}$ may be used to decompose the solution space because for any $\zeta\in\mathcal{S}$, $(\psi_j,\zeta)=0$ for ${j=1,\ldots,m}$.


Now the nonlinear terms shall be considered. The local centre manifold $W^\text{c}_{\text{loc}}$ of the equilibrium at $\mathbf{0}$ can be expressed as the sum of a linear part belonging to $\mathcal{C}$ and a nonlinear part belonging to $\mathcal{S}$,
\begin{equation}
\label{Hopf:localCM}
W^\text{c}_{\text{loc}}=\{\phi\in\mathcal{B}~|~\phi=\bm{\Phi}\mathbf{u}+\mathbf{h}(\mathbf{u})\},
\end{equation}
where $\bm{\Phi}$ is the basis given by \eqref{Hopf:basis}, $\mathbf{u}\in\mathbb{R}^m$, $\bm{\Phi}\mathbf{u}\in\mathcal{C}$, $\mathbf{h}(\mathbf{u})\in\mathcal{S}$ and $\|\mathbf{u}\|$ is sufficiently small. Thus the solutions $\mathbf{x}(t)$ to \eqref{Hopf:generalDDEtransformed} on the centre manifold satisfy $\mathbf{x}(t)=\mathbf{x}_t(0)$ where 
\begin{equation}
\label{Hopf:CMsolution}
\mathbf{x}_t(\theta)=\bm{\Phi}(\theta)\mathbf{u}(t)+\mathbf{h}(\theta,\mathbf{u}(t)).
\end{equation}
Substituting \eqref{Hopf:CMsolution} into \eqref{Hopf:operatorDE} and using \eqref{Hopf:PhiODE} and \eqref{Hopf:boundaryeqn} gives a coupled system of PDEs which must be solved for $\mathbf{u}(t)$ and $\mathbf{h}(\theta,\mathbf{u}(t))$,
\begin{align}
\label{Hopf:uhsystem}
\begin{split}
\bigg(\bm{\Phi}(\theta)&+\frac{\partial \mathbf{h}}{\partial \mathbf{u}}(\theta,\mathbf{u}(t))\bigg)\dot{\mathbf{u}}(t)\\
&=\begin{cases}
\bm{\Phi}(\theta)B\mathbf{u}(t)+\frac{\partial\mathbf{h}}{\partial\theta}(\theta,\mathbf{u}(t)),& \theta\in[-\tau,0),\\
\bm{\Phi}(0)B\mathbf{u}(t) +L\big(\mathbf{h}(\mathbf{u}(t))\big)+\mathbf{F}\big(\bm{\Phi}\mathbf{u}(t)+\mathbf{h}(\mathbf{u}(t))\big),
& \theta=0.
\end{cases}
\end{split}
\end{align}
where $L$ and $\mathbf{F}$ are defined in \eqref{L,F}, and we used the notation $\mathbf{h}(\mathbf{u}(t))$ to refer to $\mathbf{h}$ as a function in $\mathcal{S}$ for given $\mathbf{u}(t)$.
The equation for $\mathbf{u}(t)$ can now be derived using the bilinear form \eqref{Hopf:bilinearform}. Firstly, since $\mathbf{h}(\mathbf{u}(t))\in\mathcal{S}$,
\begin{equation}
\big(\bm{\Psi},\mathbf{h}(\mathbf{u}(t))\big)=0.
\end{equation}
Taking the partial derivative with respect to $\mathbf{u}$ yields
\begin{equation}
\left(\bm{\Psi},\frac{\partial\mathbf{h}}{\partial \mathbf{u}}(\mathbf{u}(t))\right)=0.
\end{equation}
Using the equations given in \eqref{L,F} and \eqref{Hopf:adjointeqns}, it can be shown that
\begin{align}
\begin{split}
\bm{\Psi}(0)&L\big(\mathbf{h}(\mathbf{u}(t))\big)
+\int_{-\tau}^{0}\bm{\Psi}(\sigma+\tau)A_1(\mu_{\rm c})\frac{\partial\mathbf{h}}{\partial\sigma}(\sigma,\mathbf{u}(t))d\sigma,\\
=&\ 
\bm{\Psi}(0) A_0(\mu_{\rm c})\mathbf{h}(0,\mathbf{u}(t))+\bm{\Psi}(\tau) A_1(\mu_{\rm c})\mathbf{h}(0,\mathbf{u}(t)) \\
&\quad
-\int_{-\tau}^{0}\bm{\Psi}'(\sigma+\tau)A_1(\mu_{\rm c})\mathbf{h}(\sigma,\mathbf{u}(t))d\sigma,\\
=&\ -B\bm{\Psi}(0)\mathbf{h}(0,\mathbf{u}(t))
-\int_{-\tau}^{0}B\bm{\Psi}(\sigma+\tau)A_1(\mu_{\rm c})\mathbf{h}(\sigma,\mathbf{u}(t))d\sigma,\\
=&\ -B\left(\bm{\Psi},\mathbf{h}(\mathbf{u}(t))\right)\\
=&\ 0.
\end{split}
\end{align}
Therefore, the bilinear form \eqref{Hopf:bilinearform} applied to $\bm{\Psi}$ and \eqref{Hopf:uhsystem} gives
\begin{equation}
\label{Hopf:ueqn}
\dot{\mathbf{u}}(t)=B\mathbf{u}(t)+\bm{\Psi}(0)\mathbf{F}\big(\bm{\Phi}\mathbf{u}(t)+\mathbf{h}(\mathbf{u}(t))\big).
\end{equation}
Substituting this into \eqref{Hopf:uhsystem} gives a system of PDEs for $\mathbf{h}(\theta,\mathbf{u}(t))$,
\begin{equation}
\label{Hopf:heqn}
\begin{aligned}
&\frac{\partial\mathbf{h}}{\partial\mathbf{u}}(\theta,\mathbf{u}(t))\Big(B\mathbf{u}(t)+\bm{\Psi}(0)\mathbf{F}\big(\bm{\Phi}\mathbf{u}(t)+\mathbf{h}(\mathbf{u}(t))\big)\Big)
+\bm{\Phi}(\theta)\bm{\Psi}(0)\mathbf{F}\big(\bm{\Phi}\mathbf{u}(t)+\mathbf{h}(\mathbf{u}(t))\big)\\
&=\begin{cases}
\frac{\partial\mathbf{h}}{\partial\theta}(\theta,\mathbf{u}(t)),&\theta\in[-\tau,0), \\
L\big(\mathbf{h}(\mathbf{u}(t))\big)+\mathbf{F}\big(\bm{\Phi}\mathbf{u}(t)+\mathbf{h}(\mathbf{u}(t))\big),
&\theta=0.
\end{cases}
\end{aligned}
\end{equation}
Standard centre manifold techniques can be used to solve \eqref{Hopf:heqn}, expanding $\mathbf{h}(\mathbf{u})$ and $\mathbf{F}$ as power series in $\mathbf{u}$ and equating coefficients. Substituting the result into \eqref{Hopf:ueqn} and expanding the right-hand side in powers of $\mathbf{u}$ gives the equation for the flow on the centre manifold. To compute the normal form for the double Hopf bifurcation, however, it is not necessary to solve the PDE, because
the coefficients of expressions up to and including cubic terms allow the normal form for the bifurcation to be identified directly using formulas available in the literature, which we now discuss.

A double Hopf bifurcation is a codimension-2 bifurcation, meaning it requires two parameters $\mu_{1,2}$ to analyse its unfolding. Let $\mu_1=\mu_{1{\rm c}}$ and $\mu_2=\mu_{2{\rm c}}$ at the bifurcation point. The normal form of the double Hopf bifurcation can be expressed in polar form in terms of two amplitudes $r_1$, $r_2$, and two phase angles $\varphi_1$, $\varphi_2$ \cite{Molnar2017} as
\begin{equation}
\label{Hopf:doublenormalform_r}
\begin{aligned}
\dot{r}_1&=b_1r_1+\left(a_{11}r_1^2+a_{12}r_2^2\right)r_1,\\
\dot{r}_2&=b_2r_2+\left(a_{21}r_1^2+a_{22}r_2^2\right)r_2,
\end{aligned}
\end{equation}
and
\begin{equation}
\label{Hopf:doublenormalform_phi}
\begin{aligned}
\dot{\varphi}_1&=\nu_1+c_{11}r_1^2+c_{12}r_2^2,\\
\dot{\varphi}_2&=\nu_2+c_{21}r_1^2+c_{22}r_2^2,
\end{aligned}
\end{equation}
where $r_{1,2} \in \mathbb{R}, \, r_{1,2}>0$ and $\varphi_{1,2} \in \mathbb{R}$. The parameters $a_{jk}$ and $c_{jk}$, $j,k\in\{1,2\}$ are known as the normal form coefficients. These can be calculated using formulae derived by \cite{Knobloch1986} from the coefficients of $\mathbf{u}$ in \eqref{Hopf:ueqn} after \eqref{Hopf:heqn} has been solved for $\mathbf{h}$ and the right-hand side has been expanded in powers of $\mathbf{u}$. 
The coefficients $b_j$, $j\in\{1,2\}$ are unfolding parameters.\footnote{Not to be confused with the parameter $b$ in the delayed HKB equation, the unfolding parameters will always have $1$ or $2$ as a subscript.} The unfolding parameters can be approximated by linear functions of the bifurcations parameters \cite{Molnar2017},
\begin{equation}
\label{Hopf:unfoldingparameters}
\begin{aligned}
b_1&=\rho_{11}\left(\mu_1-\mu_{1{\rm c}}\right)+\rho_{12}\left(\mu_2-\mu_{2{\rm c}}\right),\\
b_2&=\rho_{21}\left(\mu_1-\mu_{1{\rm c}}\right)+\rho_{22}\left(\mu_2-\mu_{2{\rm c}}\right),
\end{aligned}
\end{equation}
where 
\begin{equation}
\label{Hopf:rho}
\rho_{jk}=\text{Re}\left(\frac{\partial\lambda}{\partial\mu_k}\Biggr\rvert_{{\rm i}\nu_j}\right),
\end{equation} 
which can be calculated by differentiating the characteristic equation \eqref{Hopf:chareqn}; see Appendix~\ref{appA}. The focus of our analysis will be on the amplitudes $r_1$, $r_2$, which allow the existence and stability of the limit cycle near the bifurcation point to be investigated.

\section{Resonances}
\label{appC}

The analysis presented in this paper assumes no resonant phenomena. We address the possibility of internal resonances in the delayed HKB model \eqref{HKBeqns:full_delayed_HKB} here. The system has a $k_1:k_2$ resonance if $\frac{\nu_1}{\nu_2}=\frac{k_1}{k_2}$ for $k_1,k_2\in\mathbb{Z}^+$. If $k_1$ and $k_2$ are relatively prime integers and $k_1+k_2>4$ then the resonance is said to be weak, and no change to the normal form is required because the leading-order nonlinearities remain the same as for the non-resonant case \cite{Ma2008}. A change to the normal form is required if $k_1+k_2\leq4$.

By comparing $\nu^{(i)}$ and  $\nu^{(a)}$ at the double Hopf points, it is straightforward to show that if resonance occurs, it will be weak. The possibility of weak resonance should be investigated so that we can distinguish between resonant and non-resonant behaviour. Examples of weakly resonant and non-resonant behaviour that can arise at double Hopf bifurcations are presented by Ma et al. \cite{Ma2008}. 

 Ma et al. \cite{Ma2008} outline a method for finding the parameter values for which a $k_1:k_2$ resonant double Hopf bifurcation occurs. In the case of the delayed HKB equation, determining whether the double Hopf bifurcations found for specific parameter values correspond to resonances is not straightforward due to the complex expression for the frequencies $\nu$, given by \eqref{linear_stability_charts:nu_pm}. The relevant equations for finding the double Hopf bifurcations were solved numerically. Without analytic expressions for the parameter values at the double Hopf point, it cannot be determined with certainty that resonance is not present. However, by comparing $\nu^{(i)}$ and $\nu^{(a)}$, $k_1:k_2$ resonance can be ruled out in a finite number of cases. For example, for each of the double Hopf bifurcations found in the previous section, there is no $k_1:k_2$ resonance for $k_1,k_2\in\{1,2,\ldots,1000\}$. 

Experimentally, there is no evidence of resonance phenomena, nor are resonance phenomena mentioned by S\l{}owi\'{n}ski et al. \cite{Slowinski2016} who used the same parameter values. If experiments in the future suggest the existence of resonant phenomena for physiologically relevant parameter values, the exact parameter values which cause this can be found using the approach outlined by Ma et al. \cite{Ma2008}. However, without evidence of the phenomenon, resonance need not be discussed further here.
\section{Two-parameter continuation of invariant tori near locking}
\label{appD}

The quasi-periodic orbits have been computed using the 2D collocation algorithm in \cite{Dombovari2015,Molnar2017}, based on the method introduced in \cite{Roose2007}.
In what follows, we describe some of the trade-offs associated with this method when it is used to continue quasi-periodic orbits near locking.

The numerical continuation of quasi-periodic branches involves significantly higher errors when passing through resonances \cite{Osinga2005}.
This is the case in Figure \ref{fig:continuation} near the 1:1 strong resonance (locking) points (grey dots) where the quasi-periodic orbit collapses to the period-one limit cycle in phase quadrature (note the loss of accuracy along the red dashed lines).
We can write the governing equation (\ref{HKBeqns:full_delayed_HKB}) in the form ${\dot{\mathbf{u}} = \mathbf{f}(\mathbf{u},\mathbf{u}_\tau,\boldsymbol\mu)}$ with
${\mathbf{u} = \begin{bmatrix} x_1 & x_2 & \dot{x}_1 & \dot{x}_2 \end{bmatrix}^\top}$
and ${\mathbf{u}_\tau}:=\mathbf{u}(t-\tau)$. Computing and storing the invariant tori profiles $\mathbf{u}(\theta_1,\theta_2)$ at a given bifurcation parameter $\boldsymbol\mu$ becomes increasingly difficult as the required resolution increases. 
The numerical computation of quasi-periodic branches breaks down near the 1:1 strong resonance as it becomes increasingly difficult to satisfy the invariance relations \cite{Roose2007}
\begin{align}
\begin{split}
         \upsilon_1 \frac{\partial\mathbf{u}}{\partial\theta_1}+\upsilon_2 \frac{\partial\mathbf{u}}{\partial\theta_2}-\mathbf{f}(\mathbf{u},\mathbf{u}_\tau,\boldsymbol\mu) & = \mathbf{0},  \\
         \mathbf{u}(0,\theta_2)-\mathbf{u}(2\pi,\theta_2) & = \mathbf{0},\\
         \mathbf{u}(\theta_1,0)-\mathbf{u}(\theta_1,2\pi) & = \mathbf{0},\\
         \left\langle\frac{\partial\mathbf{u}}{\partial\theta_k},\mathbf{u}\right\rangle & = 0, \quad \quad k=1,2, \\
         \Gamma(\mathbf{u},\mathbf{u}_\tau,\upsilon_k,\boldsymbol\mu) & = 0, \quad \quad k=1,2.
\end{split}
\label{eq:invarianttori}
\end{align}
The function $\Gamma$ is given below. Pseudo-arclength \cite{Doedel2003} continuation is used to solve (\ref{eq:invarianttori}), supplemented with arclength condition. 
Equation (\ref{eq:invarianttori}) is evaluated over a Chebyshev quadrature \cite{Roose2007}.
It is not entirely clear why the breakdown of the numerical scheme arises around strong locking \cite{Osinga2005}.

Since the locking point is unreachable by continuation, approximate two-parameter $\boldsymbol\mu:=[a \,\,\tau]^\intercal$ near-locking branches - grey curves in Figure \ref{fig:continuation}(a) - were continued using the following discretization of the invariance scheme (\ref{eq:invarianttori}).
For this two-parameter continuation, the function $\Gamma$ is defined based on the fact that the numerical scheme will eventually disperse close to the 1:1 resonance, with ripples propagating higher than the $P^\textsuperscript{th}$ harmonics on the invariant torus profile given by
\begin{equation}
\begin{split}    
    \mathbf{u}(\theta_1,\theta_2) & :=\complement_{i=1,j=1}^{N,M} \mathbf{u}_{i,j}(\theta_1,\theta_2), \\
    \mathbf{u}_{i,j}(\theta_1,\theta_2) & =\sum_{p=0}^P\sum_{q=0}^P \mathbf{u}_{i,j,p,q}P_{p,q}(\epsilon_i(\theta_1),\epsilon_j(\theta_2)), \\
    P_{p,q} (\epsilon_i(\theta_1),\epsilon_j(\theta_2)) & :=P_p(\epsilon_i(\theta_1))P_q(\epsilon_j(\theta_2)),\\
    P_k(\epsilon) & :=\prod_{m=0, m\neq k}^P\frac{P\epsilon-m}{k-m}, \quad \epsilon_i(\theta_1):=\tfrac{\theta_1-\theta_{1,i}}{\Delta \theta_1}, \quad \epsilon_j(\theta_2):=\tfrac{\theta_2-\theta_{2,j}}{\Delta \theta_2},
    \end{split}    
\end{equation}
over the mesh $\theta_{1,i}:=(i-1)\Delta\theta_1$ and $\theta_{2,j}:=(j-1)\Delta\theta_2$, where $\complement$ denotes the concatenation of segmented 2D polynomial surfaces $\mathbf{u}_{i,j}$.

The different harmonics $k$ and $l$ in the spectrum are defined on the velocity profile $\dot{\mathbf{u}}_k$ of the invariant tori over both coordinates as
\begin{equation}
\begin{split}
    \mathbf{U}_{1,k}(\theta_2):=\int_0^{2\pi}\dot{\mathbf{u}}(\theta_1,\theta_2)\text{e}^{\text{i} k 2 \pi \theta_1}\text{d}\theta_1, \quad \mathbf{U}_{2,l}(\theta_1):=\int_0^{2\pi}\dot{\mathbf{u}}(\theta_1,\theta_2)\text{e}^{\text{i} l 2 \pi\theta_2}\text{d}\theta_2.
      \end{split}
\end{equation}
The main idea behind continuing near the locking curve is to lock on the relative error $\Delta$ between the main harmonics and the cumulated higher ones as
\begin{equation}
    \Gamma:=\left(\frac{\sum_{K=1,2}S_{K,1}-\sum_{K=1,2}\sum_{L=P+1}^N S_{K,L}}{\sum_{K=1,2}S_{K,1}}\right)^2-\Delta^2,
\end{equation}
with 
\begin{equation}
    S_{1,k}:=\int_0^{2\pi}\mathbf{U}_{1,k}^\textsf{H}(\theta_2) \mathbf{U}_{1,k}(\theta_2)\text{d}\theta_2, \quad S_{2,l}:=\int_0^{2\pi}\mathbf{U}_{2,l}^\textsf{H}(\theta_1) \mathbf{U}_{2,l}(\theta_1)\text{d}\theta_1.
\end{equation}
Corresponding derivatives can be computed using the chain rule and Newton-Raphson correction can be performed quite effectively.

\bibliography{MScSources}
\bibliographystyle{plain}

 \end{document}